\documentclass{elsarticle}

\usepackage[utf8]{inputenc}
\usepackage{amsmath}
\usepackage{amssymb, amsthm, units, bm, esint}
\usepackage[T1]{fontenc}
\usepackage{srcltx}
\usepackage{graphicx}
\usepackage{subfig}
\usepackage{accents}
\usepackage{natbib}
\usepackage{hyperref}
\usepackage[colorinlistoftodos]{todonotes}

\newtheorem{thm}{Theorem}[section]
\newtheorem*{thm*}{Theorem}

 \newtheorem*{lem*}{Lemma}
 \newtheorem{prop}[thm]{Proposition}
 \newtheorem{prop*}{Proposition}[section]
 \theoremstyle{definition}
 
 \theoremstyle{remark}

 \numberwithin{equation}{section}

\def\tens#1{\mathbb{#1}}
\def\vec#1{\boldsymbol{#1}}
\def \R{\mathbb{R}}
\def \N{\mathbb{N}}

\def\bD{\tens{D}}
\def\bG{\tens{G}}
\def\bO{\tens{O}}

\def\bI{\tens{I}}

\def\bF{\tens{F}}

\def\bB{\tens{B}}
\def\bA{\tens{A}}
\def\bC{\tens{C}}
\def\bH{\tens{H}}

\def\bU{\tens{U}}
\def\bV{\tens{V}}
\def\bQ{\tens{Q}}

\def\bu{\vec{u}}
\def\bv{\vec{v}}
\def\bw{\vec{w}}

\def\bn{\vec{n}}
\def\bb{\vec{b}}
\def\bz{\vec{z}}

\def\be{\begin{equation}}
\def\ba{\begin{array}}
\def\ea{\end{array}}
\def\ee{\end{equation}}

\DeclareMathOperator{\diver}{div}
\DeclareMathOperator{\Diver}{Div}
\DeclareMathOperator{\tr}{tr}

\DeclareMathOperator{\adj}{adj}

\begin{document}

\begin{frontmatter}

\title{On three-dimensional flows of viscoelastic fluids of Giesekus type\tnoteref{mytitlenote}}
\tnotetext[mytitlenote]{Miroslav Bul\'{i}\v{c}ek,  Tom\'{a}\v{s} Los and Josef M\'{a}lek thank to the Project 20-11027X financed by the Czech Science Foundation. Miroslav Bul\'{\i}\v{c}ek and Josef M\'{a}lek are members of the Ne\v{c}as Center for Mathematical Modelling.}

\author[Bulicek-address]{Miroslav Bul\'{i}\v{c}ek}
\ead{mbul8060@karlin.mff.cuni.cz}
\author[Bulicek-address]{Tom\'{a}\v{s} Los}
\ead{los@karlin.mff.cuni.cz}
\author[Bulicek-address]{Josef M\'{a}lek}
\ead{malek@karlin.mff.cuni.cz}

\address[Bulicek-address]{Charles University, Faculty of Mathematics and Physics,  Sokolovsk\'{a}~83, 186~75 Praha~8, Czech Republic}


\begin{abstract}
Viscoelastic rate-type fluids are popular models of choice in many applications involving flows of fluid-like materials with complex micro-structure. A  well-developed mathematical theory for the most of these classical fluid models is however missing. The main purpose of this study is to provide a complete proof of long-time and large-data existence of weak solutions to unsteady internal \emph{three-dimensional} flows of Giesekus fluids subject to a no-slip boundary condition. As a new auxiliary tool, we provide the identification of certain biting limits in the parabolic setting,  presented here within the framework of evolutionary Stokes problems. We also generalize the long-time and large-data existence result to higher dimensions, to viscoelastic models with multiple relaxation mechanisms and to viscoelastic models with different type of dissipation. 
\end{abstract}


\end{frontmatter}

\section{Introduction}

It has been experimentally observed and numerically tested that viscoelastic rate-type fluid models are capable of capturing the properties of materials with complicated micro-structure (polymers, bio-materials, geo-materials) and can thus serve as a mathematical setting suitable for predicting flows of such materials. To improve the efficiency and robustness of (finite-element based) numerical schemes designed to compute flows governed by these models, it is important to have at one's disposal a well-founded body of mathematical theory for such models. For the most of the classical viscoelastic rate-type fluid models such a theory is however missing. Here we aim to contribute to filling in this lacuna by focusing on three-dimensional unsteady flows.  

More specifically, we present a proof of the long-time and large-data existence of weak solutions to the Giesekus model in three spatial dimensions. 
We thus extend and significantly strengthen the results established in the recent study~\cite{BuMaLo22} focused on two-dimensional/planar flows. More precisely, in \cite{BuMaLo22}, following ideas from~\cite{Mas} where N. Masmoudi outlined the proof of weak sequential stability of hypothetical weak solutions to the Giesekus viscoelastic rate-type fluid model (see \cite{GIESEKUS}), a complete proof of the long-time and large-data existence of weak solutions to the Giesekus model in two spatial dimensions is given. In addition, also in \cite{BuMaLo22}, this proof is extended to a mixture of two Giesekus models (with two different relaxation times); such a model can be viewed as a generalization of the classical Burgers model (which is the most popular model in the class of viscoelastic rate-type fluid models of the second order). We also documented the importance of these models by recalling their recent applications in modelling responses of such diverse materials as rubber, asphalt binder or vitreous gel, see \cite{Reh, Nar, derivold, Sk, Tuma18}.  Finally, following \cite{MRT}, we also put these viscoelastic rate-type models of higher order into a unified thermodynamic setting. We refer the reader to \cite{BuMaLo22} for more details. 

\medskip

We shall now formulate the main results of this paper. 

\smallskip

\noindent Throughout the whole text (unless explicitly stated otherwise), we assume that 
\begin{equation}
    \label{p000}
    \textrm{\begin{minipage}{10cm}
    $\bullet$ $\Omega\subset \R^d$ is a  bounded open connected set with a Lipschitz boundary $\partial \Omega$ and with the dimension $d\ge 2$, \\ 
    $\bullet$ $T\in(0,\infty)$ is the length of an arbitrary but fixed time interval, \\ $\bullet$ $Q_T:=(0,T)\times~\Omega$.\end{minipage}}
\end{equation}

\noindent 
Given the parameters $G,\tau, \rho, \nu\in (0,\infty)$ and the external body forces $\vec{f}:Q_T\to \R^d$, we look for functions $(\bv, p, \bF, \bB): Q_T\to \R^d\times \R \times \R^{d\times d} \times \R^{d\times d}$ satisfying the following governing equations and the natural constraints in $Q_T$:
\begin{align}
\label{1Burg}
\diver \bv&=0,\\
\rho\partial_t \bv + \rho\diver (\bv\otimes \bv) + \nabla p - \nu\diver (\nabla \bv+(\nabla \bv)^T)-G \diver\bB&=\rho\vec{f},\\
\label{vB3lambda}
\partial_t \bB +\mbox{Div} (\bB\otimes \bv) - \nabla \bv \bB -\bB (\nabla \bv)^T+ \frac{1}{\tau}\left(\bB^2 - \bB\right)&=\tens{O},\\
\bF \bF^T&=\bB, \\
\det \bF&>0,  \label{detFform}
\end{align}
together with the homogeneous Dirichlet boundary condition for $\bv$:
\begin{equation}\label{Burgbc}
\bv=\vec{0} \quad \mbox{on } (0,T)\times \partial \Omega,
\end{equation}
and the initial conditions for $\bv$, $\bF$ and $\bB$:
\begin{equation}\label{ic}
\bv(0,\vec{\cdot})=\bv_0, \quad \bF(0,\vec{\cdot})=\bF_{0}, \quad \bB(0,\vec{\cdot})=\bB_{0}:= \bF_{0} \bF_{0}^T \quad \mbox{in } \Omega.
\end{equation}
The problem \eqref{1Burg}--\eqref{ic} is compatible with the laws of continuum thermodynamics and satisfies, for all $t>0$, the following energy identity (see \cite[Sect.~2]{BuMaLo22}): 
\begin{equation} 
\label{pepa98}
\begin{split}
\int_{\Omega} \rho \frac{|\bv(t)|^2}{2} &+ \frac{G}{2}( \tr \bB (t) -d-\ln \det \bB(t)) +  \int_0^t \int_{\Omega}  2\nu |\bD|^2 + \frac{G}{2\tau} |\bB -\bI|^2  \\
&= \int_{\Omega} \rho\frac{|\bv_0|^2}{2} + \frac{G}{2} (\tr \bB_{0} - d- \ln \det \bB_{0}) + \int_0^t \int_{\Omega} \rho \vec{f}\cdot \bv.
\end{split}
\end{equation}
Here, $\bv$, $p$, $\bB$ and $\vec{f}$ are smooth enough, $2\bD:=\nabla \bv +(\nabla \bv)^T$, $\bI$ is the identity tensor and $\tr \bB = |\bF|^2$.

\smallskip\noindent 
Our first main result can be informally formulated as follows.
\begin{thm*}
For an arbitrary $\Omega$ and $T$ satisfying  \eqref{p000} and for arbitrary data $\bv_0$, $\bF_0$ and $\vec{f}$ that enables the right-hand side of \eqref{pepa98} to be finite, there exists a global-in-time weak solution to the problem \eqref{1Burg}--\eqref{ic}.
\end{thm*}
\noindent A precise formulation of this result  that in particular concerns \emph{three-dimensional} flows is given in Theorem~\ref{Burgmain} below.

\medskip

Although there are similarities with the method outlined in~\cite{Mas} and with the two-dimensional case treated in a full detail rigorously in~\cite{BuMaLo22}, there is one critical point, which makes the analysis in the three-dimensional setting much more challenging and difficult. Namely, one cannot use the energy \emph{equality} when proving the compactness of $\bF$ in $(L^2(Q_T))^{d\times d}$ and one must find another way to identify the limit in the nonlinear terms. This leads us to developing a new approach that can be formulated as a result \emph{per se}. This novel method then forms the second main new result of this study, see Theorem~\ref{bitingT} below. 

The structure of this study is the following. In the remaining part of this introductory section we recall the basic notation used in the paper. We also summarize several fundamental tools needed in the proofs of the main theorems. The introductory section ends with rigorous formulations of the main results in Theorem~\ref{Burgmain} (existence of a large-data global-in-time weak solution for the Giesekus model) and in Theorem~\ref{bitingT} (identification of certain biting limits). Section~\ref{Proof1} is devoted to the proof of Theorem~\ref{bitingT} and in Section~\ref{Proof2} we provide the complete proof of Theorem~\ref{Burgmain}. Various  generalizations of the result (to higher dimensions, to viscoelastic models with multiple relaxation mechanisms and to viscoelastic models with different type of dissipation) are formulated in Section~\ref{Extensions}.

\subsection{Basic notation}

The symbol $``\cdot"$ denotes the scalar product of two vectors, while the symbol $``:"$ is used to denote the scalar product of two tensors. The symbol $``~\otimes~"$ denotes the tensor product of two vectors, i.e. $(\bb\otimes\bz)_{ij}:= b_i z_j$. For a matrix $\bA=\{A_{ij}\}_{i,j=1}^d$ and a vector $\bb=(b_1,\ldots,b_d)$ we define the third order tensor $\bA\otimes ~\bb= \left\{(\bA\otimes \bb)_{i j k}\right\}_{i,j,k=1}^d$ as $(\bA\otimes \bb)_{ijk}:= A_{ij} b_k$. For a matrix-valued function $\bA=(A_{ij})_{i,j=1}^d$ and a vector-valued function $\bb=(b_1,...,b_d)$ we define  the operator $\Diver$ acting on the third order tensor $\bA\otimes \bb$ as
$$
\Diver (\bA\otimes \bb):= \sum_{j=1}^d \partial_{x_j} (b_j \bA).
$$

As stated above, we assume that $\Omega\subset \R^d$ and $T\in(0,\infty)$ satisfy \eqref{p000}. The symbol $|\cdot|$, as usual, has several different meanings: it denotes the Euclidean norm of a vector, the Frobenius norm of a tensor, and also the Lebesgue measure of a given measurable subset of $\R^n$, where $n=d$ or $n=d+1$. Further, if $O\subset \R^n$  is a measurable set, then $\chi_O$ denotes the characteristic function of $O$.

For $q\in [1, \infty]$, the symbols $\|\cdot\|_q$ and $\|\cdot\|_{m,q}$, $m\in \N$, stand for the norms in the usual Lebesgue and Sobolev spaces $L^q (\Omega)$ and $W^{m,q} (\Omega)$ (or their multi-dimensional variants). Next, for $O\subset \R^n$ open, $\mathcal{C}_c^\infty(O)$ denotes the space of smooth functions compactly supported in~$O$. Finally, we set
\begin{align*}
L_{\bn, \diver}^q &:=  \overline{\{\bu \in (\mathcal{C}_c^\infty (\Omega))^d; \diver \bu = 0 \mbox{ in }  \Omega\}}^{\|\cdot\|_q},\\
W_{0}^{1,q} (\Omega)&:=  \{u\in W^{1,q} (\Omega); u =0 \mbox{ on } \partial \Omega\},\\
W_{\vec{0},\diver}^{1,q} &:= \{\bu\in (W^{1,q} (\Omega))^d; \bu=\vec{0} \mbox{ on }  \partial \Omega; \diver \bu = 0  \mbox{ in } \Omega\},
\end{align*}
and we equip these spaces with the norms (thanks to Poincar\'{e}'s inequality)
$$
\|u\|_{L^q_{\bn,\diver}}:= \|u\|_q,\quad
\|u\|_{W_0^{1,q}(\Omega)}:= \|\nabla u\|_q,\quad \|\bu\|_{W_{\vec{0},\diver}^{1,q}}:= \|\nabla \bu\|_q.
$$
The symbol $\mathcal{M}(\overline{Q_T})$ stands for the space of Radon measures defined on the closure of $Q_T$. More details regarding rather standard notation and conventions are given in \cite{BuMaLo22}. The standard abbreviations \emph{a.e.} and \emph{a.a.} stand for \emph{almost everywhere} and \emph{almost all}. 

\subsection{Mathematical and convergence tools}\label{MCTS2}
Here, we provide a list of results that are used to prove the main theorems of the paper. First, we recall the parabolic Lipschitz approximation method, which allows us to identify certain biting limits.
\begin{prop}[Parabolic Lipschitz truncation, {\cite[Theorem 3.1]{BuBuSch19}}]
\label{BuLip}
Let $T\in (0,\infty)$ and $\tilde{\Omega}\subset \overline{\tilde{\Omega}}\subset\Omega$ satisfy the same assumptions as $\Omega$, as specified in 
\eqref{p000}. Assume that $\{\bw_k\}_{k\in \N}\subset L^2(0,T;(W^{1,2}(\Omega))^d)$ and $\{\bG_k\}_{k\in \N}\subset (L^2(Q_T))^{d\times d}$ are such that
\begin{equation}\label{BuLipasest}
   \sup_{k\in \N} {} \left(\|\bG_k\|_{(L^2(Q_T))^{d\times d}}+\|\bw_k\|_{L^2(0,T;(W^{1,2}(\Omega))^d)}\right)< \infty
\end{equation}
and the following identity holds for any $k\in \N$ and all $\vec{\varphi}\in (\mathcal{C}_c^\infty((-\infty, T)\times \Omega))^d$:
\begin{equation}\label{BuLipaseq}
    \int_{Q_T}\bw_k\cdot \partial_t \vec{\varphi}-\int_{Q_T} \bG_k:\nabla \vec{\varphi} =0.
\end{equation}
Then, for any $\Lambda> 1$, there exist a sequence  $\{\bw_{\Lambda,k}\}_{k\in \N}\subset L^2(0,T;(W^{1,2}(\tilde{\Omega}))^d)$ and a constant $C\in (0,\infty)$ independent of $k$ and $\Lambda$ such that, for all $k\in~\N$, the following holds:
\begin{align}
        \label{estnablwlambdak}&\|\nabla \bw_{\Lambda,k}\|_{(L^\infty((0,T)\times\tilde{\Omega}))^{d\times d}}+\|\partial_t \bw_{\Lambda,k}\cdot (\bw_{\Lambda,k}-\bw_ k)\|_{L^2((0,T)\times\tilde{\Omega})}\leq C\Lambda^{2^{\Lambda}},\\
     \label{estdtwlambdak}&\|\nabla \bw_{\Lambda,k}\|_{(L^2((0,T)\times\tilde{\Omega}))^{d\times d}}+\sqrt{\Lambda}\|\partial_t \bw_{\Lambda,k}\cdot (\bw_{\Lambda,k}-\bw_ k)\|_{L^1((0,T)\times\tilde{\Omega})}\le C,\\
     \label{estwlambdakk}&\|\bw_{\Lambda, k}\|_{(L^s((0,T)\times\tilde{\Omega}))^d}\leq C \|\bw_k\|_{(L^s((0,T)\times\tilde{\Omega}))^d} \qquad \textrm{ for all } s \in [1,\infty),
\end{align}
and,  for all $\psi\in \mathcal{C}^{0,1}_c((0,T)\times\tilde{\Omega})$, the following identity holds:
\begin{equation}\label{BuLipeq}
\begin{split}
    \int_{Q_T} \bG_k: \nabla (\bw_{\Lambda,k}\psi)= &-\frac{1}{2}\int_{Q_T} (|\bw_{\Lambda,k}|^2-2\bw_k\cdot \bw_{\Lambda,k})\partial_t\psi\\
    &-\int_{Q_T}\partial_t \bw_{\Lambda, k}\cdot (\bw_{\Lambda,k}-\bw_k)\psi.
\end{split}
\end{equation}
In addition, defining the set $\mathcal{O}_{\Lambda, k}$ as
\begin{equation}\label{defOLambdak}\mathcal{O}_{\Lambda, k}:= \{[t,\vec{x}]\in (0,T)\times\tilde{\Omega};\ \bw_{\Lambda,k}(t,\vec{x})\neq \bw_k(t,\vec{x})\},\end{equation}
we get the estimate
\begin{align}
\label{estOLambdak}
|\mathcal{O}_{\Lambda, k}|&\leq \frac{C}{\Lambda}.
\end{align}
\end{prop}

Next, since our method is based on the use of the biting limit, we recall the fundamental property of sequences bounded in $L^1$.
\begin{prop}[Chacon's biting lemma, see \cite{BaMu89}]\label{biting}
Let $Q\subset \mathbb{R}^{d+1}$ be a measurable set. Assume that a~sequence $\{f_k\}_{k\in \N}$ satisfies
\begin{equation}\label{L1bound}
\sup_{k\in \N} \|f_k\|_{L^1(Q)}<\infty.
\end{equation}
Then there exist: \emph{(i)} $f \in L^1(Q)$, \emph{(ii)} a non-decreasing family of measurable sets $\{E_j\}_{j\in \N}$, $E_j\subset Q$ for all $j\in \N$, satisfying $\lim_{j\to\infty}|Q\setminus E_j|=0$,
and \emph{(iii)} a subsequence of $\{f_k\}_{k\in \N}$, which we do not relabel, 
such that, for all $j\in \N$, the following holds:
\begin{equation}\label{bitingesta}
    f_k \rightharpoonup f \quad \textrm{weakly in } L^1(E_j).
\end{equation}
\end{prop}

In what follows, whenever we write
$$
f_k \rightharpoonup f \textrm{ in the biting sense in }L^1(Q),
$$
this means that, for the sequence $\{f_k\}_{k\in \N}$ fulfilling \eqref{L1bound}, we are applying Chacon's biting lemma (Proposition~\ref{biting}) and find the corresponding $f\in L^1(Q)$ and a sequence of sets $\{E_j\}_{j\in \N}$ with the properties described above (and a suitable subsequence of $\{f_k\}_{k\in \N}$ that is not relabeled) such that \eqref{bitingesta} holds. Moreover, it then also follows from the characterization of $L^1$-weakly converging sequences that $\{f_k\}_{k\in \N}$ is uniformly equi-integrable on $E_j$ for any~$j\in \N$. This means that, for any~$\varepsilon >0$ and any~$j\in \N$, there exists a $\delta>0$ such that, for all measurable sets $U\subset E_j$ fulfilling $|U|\le \delta$, there holds
\begin{equation}
\label{unifL1}
\sup_{k\in \N} \int_{U}|f_k| \le \varepsilon.
\end{equation}
To be precise, the uniform bound~\eqref{L1bound} implies, for a suitable subsequence of $\{f_k\}_{k\in \N}$ that we do not relabel, the validity of the estimate \eqref{unifL1}, which, for a further not relabeled subsequence of $\{f_k\}_{k\in \N}$, implies the convergence result~\eqref{bitingesta}. On the other hand, the convergence result~\eqref{bitingesta} implies the validity of \eqref{unifL1}.

\medskip

The next two results are based on the classical theory for the Stokes system.
\begin{prop}[Stokes problem I, see \cite{Wolf} or {\cite[Appendix A.1.]{BuMaLo22}}]\label{recpres1}
Let \eqref{p000} \\hold and let the sequences $\{\bH_k\}_{k\in \N}\subset (L^2(Q_T))^{d\times d}$ and  $\{\bu_{0,k}\}_{k\in \N}\subset L^2_{\bn,\diver}$ be given. Then, for every fixed $k\in \N$, there exists a unique  $\bu_k\in L^2(0,T;W_{\vec{0},\diver}^{1,2})$, such that $\partial_t \bu_k\in L^2(0,T;(W_{\vec{0},\diver}^{1,2})^*)$, and, for all $\bw\in W_{\vec{0},\diver}^{1,2}$ and a.a. $t\in (0,T)$, the following equation holds (here $\bD_k:=\frac{\nabla \bu_k+(\nabla \bu_k)^T}{2}$)
 \begin{equation}
     \langle \partial_t \bu_k, \bw\rangle+\int_{\Omega}\bD_k:\nabla \bw=-\int_{\Omega}\bH_k:\nabla \bw\label{dva12}
 \end{equation}
with the initial condition $\bu_{0,k}$ fulfilled in the sense
\begin{equation}
     \lim_{t\to0+} \|\bu_k(t)-\bu_{0,k}\|_2=0.\label{dva13}
\end{equation}
Also, for all $k\in \N$, there exists a unique (up to a constant) pressure $p_k:=\partial_t p_{h,k}+\tilde{p}_k$ such that, for every smooth $\tilde{\Omega}$ satisfying 
$\tilde{\Omega}\subset\overline{\tilde{\Omega}}\subset \Omega$, there holds
 \begin{align}
     p_{h,k} &\in L^2(0,T;W^{2,2}(\tilde{\Omega})),\quad \Delta p_{h,k}=0\mbox{ a.e. in }(0,T)\times\tilde{\Omega},\\
     \tilde{p}_k&\in L^2((0,T)\times\tilde{\Omega}),\\
     \partial_t (\bu_k+\nabla p_{h,k})&\in L^2(0,T;((W_{0}^{1,2}(\tilde{\Omega}))^d)^*),
 \end{align}
 and, for all $\bw\in (W^{1,2}_0(\tilde{\Omega}))^d$  and a.a. $t\in (0,T)$, we have
 \begin{equation}
     \langle\partial_t (\bu_k+\nabla p_{h,k}),\bw\rangle-\int_{\Omega}\tilde{p}_k\diver \bw+\int_{\Omega} \bD_k:\nabla \bw=-\int_{\Omega}\bH_k:\nabla \bw.
 \end{equation}
Moreover, if $\sup_{k\in \N}\|\bH_k\|_{(L^2(Q_T))^{d\times d}}<\infty$ and $\sup_{k\in \N} \|\bu_{0,k}\|_{L^2_{\bn,\diver}}<\infty$, then there exist $\bu$, $\bu_0$, $p_h$, $\tilde{p}$ and $\bH$ such that, for suitable subsequences, which we do not relabel, there holds
\begin{align}
     \bH_k&\rightharpoonup \bH && \mbox{weakly in } (L^2(Q_T))^{d\times d},\\
     \label{ukweak1}\bu_k&\rightharpoonup \bu && \mbox{weakly in } L^2(0,T; W_{\vec{0},\diver}^{1,2}),\\
     \label{dtukweak}\partial_t \bu_k&\rightharpoonup\partial_t \bu&& \mbox{weakly in } L^2(0,T;(W_{\vec{0},\diver}^{1,2})^*),\\
     \bu_k&\to \bu&& \mbox{strongly in } (L^2(Q_T))^d,\\
     \bu_{0,k}&\rightharpoonup \bu_0&& \mbox{weakly in } L^2_{\bn,\diver},\\
    \label{phkstrong} p_{h,k}&\to p_{h} && \mbox{strongly in } L^2(0,T;W^{2,2}(\tilde{\Omega})),\\
    \label{tildepweak} \tilde{p}_{k}& \rightharpoonup \tilde{p} && \mbox{weakly in } L^2((0,T)\times \tilde{\Omega}).
\end{align}
In addition, we have that $\Delta p_h=0$ a.e. in $(0,T)\times \tilde{\Omega}$,
\begin{equation}\label{contukphk}
    \bu_{k}, \bu \in \mathcal{C}([0,T]; L^2_{\bn,\diver}), \quad \nabla p_{h,k}, \nabla p_{h} \in \mathcal{C}([0,T]; (L^2(\tilde{\Omega}))^d)
\end{equation}
and
\begin{equation}\label{icukphk}
    \bu(0,\vec{\cdot})=\bu_0\ \mbox{ a.e. in } \Omega.
\end{equation}
\end{prop}


\begin{prop}[Stokes problem II]\label{recpres2}
Let \eqref{p000} hold. Assume that the sequences $\bH_k\in (L^2(Q_T))^{d\times d}$ and $\bu_k  \in L^2(0,T;W_{\vec{0},\diver}^{1,2})$ 
are such that, for all $\bw\in W_{\vec{0},\diver}^{1,2}$ and a.a. $t\in(0,T)$, the following identity holds (here $\bD_k:=\frac{\nabla \bu_k+(\nabla \bu_k)^T}{2}$):
\begin{equation}\label{SII}
     \langle \partial_t \bu_k, \bw\rangle+\int_{\Omega} \bD_k: \nabla \bw=-\int_{\Omega} \bH_k:\nabla \bw
\end{equation}
and $\bu_k(0)=\vec{0}$ a.e. in $\Omega$. Furthermore, let, for some $q>1$, 
$$
\begin{aligned}
\nabla \bu_k &\rightharpoonup \nabla \bu &&\textrm{weakly in }(L^2(Q_T))^{d\times d},\\
\bH_k &\to \bH &&\textrm{strongly in } (L^q(Q_T))^{d\times d}.
\end{aligned}
$$
Then, for any $r\in [1,2)$,
\begin{align}\label{Bi21}
\nabla \bu_k &\to \nabla \bu &&\textrm{strongly in }(L^r(Q_T))^{d\times d}.
\end{align}
\end{prop}
The above proposition is a direct consequence of results obtained in \cite{DiRuWo10}, see also \cite{BlMaRa20,BuGwMaSw12} or \cite{KSol}.
%
%
The last result, we recall here, is a kind of the Gronwall-like lemma for certain differential inequalities.
\begin{prop}\label{convae}
Let \eqref{p000} hold and $\bu\in L^2(0,T; W_{\vec{0},\diver}^{1,2})$ and $L\in L^2(Q_T)$ be given. Assume that, for a non-negative $f\in L^2(Q_T)$, the following holds:  
\begin{equation}
\begin{split}
     -\int_{Q_T} f \partial_t \varphi-\int_{Q_T}f \bu\cdot \nabla \varphi &\leq \int_{Q_T} L f \varphi \\
     &\textrm{for all non-negative } \varphi\in \mathcal{C}_c^\infty((-\infty,T)\times \Omega).
\end{split}
\end{equation}
Then, $f=0$ a.e. in $Q_T$.
\end{prop}
\begin{proof}
The ideas used in the proof of the above proposition can be traced back to the classical result of R.~DiPerna and P.-L.~Lions in~\cite{DiPerna1989}. However, the proof presented in~\cite{DiPerna1989} works  only if $L\in L^1(0,T; L^{\infty}(\Omega))$, which is not the case in Proposition~\ref{convae}. Nevertheless,  statements very similar to 
Proposition~\ref{convae} are known and used in the theory of fluid flows thanks to P.-L.~Lions, see e.g. \cite{Li} and E.~Feireisl, see e.g. \cite{Fe04}. However, we were not able to find a proper reference describing the case treated in Proposition~\ref{convae}. This is why we refer to \cite[Step 3: the proof of the compactness of $\{\bF_\varepsilon\}$ in $(L^2(Q_T))^{2\times2}$ and Remark below Step 3]{BuMaLo22} for its proof.
\end{proof}

\subsection{Formulation of the Main Result}
\label{FormS3}

The main result of this paper concerns the global-in-time and large-data existence of weak solutions to the problem \eqref{1Burg}--\eqref{ic} and it is formulated in the following theorem. Let us note that, throughout the paper, we use the notation $\bD:=\frac{\nabla \bv+(\nabla \bv)^T}{2}.$ To simplify the presentation, we formulate and prove the existence result only in dimension $d=3$. However, note that the result and also the proof can be straightforwardly adapted also to a general $d$-dimensional setting, the only difference being that one then uses a different interpolation inequality for the velocity field $\bv$, which consequently results in different function spaces for $\partial_t\bv$ and $\partial_t \bB$. See also Sect.~\ref{Extensions}.
\begin{thm}
\label{Burgmain}
Let $\Omega\subset \R^3$ and $T$ satisfy \eqref{p000} and let $G, \tau, \rho, \nu\in(0,\infty)$ be given constants. Let $\vec{f}\in L^2(0,T;(W_{\vec{0},\diver}^{1,2})^*)$, $\bv_0\in L^2_{\bn, \diver}$, $\bF_{0}\in (L^2(\Omega))^{3\times 3}$, $\det \bF_{0}>0$ a.e. in $\Omega$ and $\ln \det \bF_{0}\in L^1(\Omega)$.  Then there exists a weak solution to \eqref{1Burg}--\eqref{ic}, i.e. there exists a triple $(\bv, \bF,  \bB)$ satisfying
\begin{align*}
\bv & \in \mathcal{C}_{weak}([0,T]; L^2_{\bn, \diver}) \cap L^2(0,T; W_{\vec{0},\diver}^{1,2}), \\
\partial_t \bv &\in L^{\frac{4}{3}}(0, T; (W_{\vec{0},\diver}^{1,2})^*),\\
\bB & \in \mathcal{C}([0,T]; (L^1(\Omega))^{3\times 3}) \cap (L^2(Q_T))^{3\times3}, \\
\partial_t \bB & \in L^1 (0,T; ((W^{1,4}(\Omega))^{3\times 3})^*),\\
\bF \bF^T&= \bB \mbox{ a.e. in } Q_T,\\
\det \bF &> 0 \mbox{ a.e. in } Q_T,
\end{align*}
such that, for all $\bw \in W_{\vec{0},\diver}^{1,2}$, all $\bA \in (W^{1,4}(\Omega))^{3\times 3}$ and a.a. $t\in (0,T)$, the following identities hold:
\begin{align}
\label{vBweaklambda1}
\begin{aligned}
\rho\langle \partial_t \bv, \bw \rangle &- \rho\int_{\Omega} (\bv \otimes \bv): \nabla \bw + \int_{\Omega} (2\nu\bD+G\bB): \nabla \bw\\ & -\langle \vec{f}, \bw\rangle = 0,
\end{aligned} \\
\label{vBweaklambda2}
\begin{aligned}
    \langle \partial_t \bB, \bA \rangle &-\int_{\Omega} (\bB \otimes \bv): \nabla \bA - \int_{\Omega} (\nabla \bv \bB+\bB (\nabla \bv)^T): \bA \\
    &+ \frac{1}{\tau}\int_{\Omega} (\bB^{2} -  \bB): \bA = 0.
\end{aligned} 
\end{align}
The initial conditions $\bv_0$, $\bF_{0}$, $\bB_{0}:= \bF_{0} \bF_{0}^T$ are attained as follows:
\begin{equation}
\label{invFBweak}
\lim_{t\to 0+} \left(\|\bv(t)-\bv_0\|_2 +\|\bF(t)-\bF_{0}\|_2 + \|\bB(t)-\bB_{0}\|_1\right)= 0.
\end{equation}
\end{thm}

The second main result is a new convergence tool for the Stokes (or parabolic) systems, identifying the biting limit of the ``energy"/``dissipation" functional.
\begin{thm}\label{bitingT}
Let $\Omega$ and $T$ satisfy \eqref{p000} and  $\{\tilde{\bw}_k\}_{k\in \N}$ and $\{\tilde{\bG}_k\}_{k\in \N}$ be sequences satisfying
\begin{align}\label{bitingC}
\begin{aligned}
\tilde{\bw}_k &\rightharpoonup \tilde{\bw} &&\textrm{weakly in }L^2(0,T;W_{\vec{0},\diver}^{1,2}),\\
\partial_t \tilde{\bw}_k &\rightharpoonup \partial_t\tilde{\bw} &&\textrm{weakly in }L^2(0,T;(W_{\vec{0},\diver}^{1,2})^*),\\
\tilde{\bG}_k &\rightharpoonup \tilde{\bG} &&\textrm{weakly in }(L^{2}(Q_T))^{d\times d}
\end{aligned}
\end{align}
and fulfilling, for a.a. $t\in (0,T)$ and all $\varphi \in W_{\vec{0},\diver}^{1,2}$, the following identity
\begin{equation}\label{formB}
\langle \partial_t \tilde{\bw}_k, \vec{\varphi}\rangle + \int_{\Omega} \tilde{\bG}_k : \nabla \vec{\varphi} =0.
\end{equation}
Then there holds
\begin{align}\label{Br}
\tilde{\bG}_k : \nabla \tilde{\bw}_k \rightharpoonup \tilde{\bG} : \nabla \tilde{\bw} \qquad \textrm{in the biting sense in } L^1(Q_T).
\end{align}
\end{thm}
The main novelty of the above theorem is  the fact that \eqref{Br} holds for the biting limit. Note that it is not difficult to prove via a standard energy method that
\begin{align*}
\tilde{\bG}_k : \nabla \tilde{\bw}_k \rightharpoonup^* \tilde{\bG} : \nabla \tilde{\bw} \qquad \textrm{weakly$^*$ in } \mathcal{M}(Q_T),
\end{align*}
i.e., $\tilde{\bG}_k : \nabla \tilde{\bw}_k$ converges to $\tilde{\bG} : \nabla \tilde{\bw}$ weakly$^*$ in the space of the Radon measures. However, the biting limit and the limit in measures may in general differ significantly (on a set of positive measure). Thus, the main novelty here is the observation that \emph{the limit in measures coincides with the limit in the biting sense}.

\section{Proof of Theorem~\ref{bitingT}}\label{Proof1}
First, we define the new quantities
\begin{equation}\label{Berlin}
\begin{aligned}
\bu_k&:=\tilde{\bw}_k - \tilde{\bw}, \qquad \bu_{0,k}:=\tilde{\bw}_k(0)-\tilde{\bw}(0),\\
\bD_k&:=\frac12 \left(\nabla \bu_k + (\nabla \bu_k)^T\right),\\
\bH_k&:= \tilde{\bG}_k-\tilde{\bG} - \bD_k.
\end{aligned}
\end{equation}
Then, it follows from the assumptions of Theorem~\ref{bitingT} that the quadruple $(\bu_k, \bu_{0,k}, \bD_k, \bH_k)$ satisfies \eqref{dva12}--\eqref{dva13}. Note that $\bu_k$ is uniquely defined from $\bH_k$ and $\bu_{0,k}$  by using Proposition~\ref{recpres1} and, in addition, we have the following convergence results as $k\to \infty$:
\begin{align}
     \label{Hkweak}\bH_k&\rightharpoonup \tens{O} && \mbox{weakly in } (L^2(Q_T))^{d\times d},\\
     \label{ukweak2}\bu_k&\rightharpoonup \vec{0} && \mbox{weakly in } L^2(0,T; W_{\vec{0},\diver}^{1,2}),\\
     \partial_t \bu_k&\rightharpoonup \vec{0}&& \mbox{weakly in } L^2(0,T;(W_{\vec{0},\diver}^{1,2})^*),\\
     \label{ukstrong2}\bu_k&\to \vec{0}&& \mbox{strongly in } (L^2(Q_T))^d,\\
     \bu_{0,k}&\rightharpoonup \vec{0}&& \mbox{weakly in } (L^2(\Omega))^d.
\end{align}
In addition, still by Proposition~\ref{recpres1}, we can find the corresponding pressures $\partial_t p_{h,k}$ and $\tilde{p}_k$ such that, for every  smooth $\hat{\Omega}\subset\overline{\hat{\Omega}}\subset \Omega$, there holds
\begin{align}
p_{h,k} &\in L^2(0,T;W^{2,2}(\hat{\Omega})),\quad \Delta p_{h,k}=0\mbox{ a.e. in }(0,T)\times\hat{\Omega},\\
\tilde{p}_k&\in L^2((0,T)\times\hat{\Omega}),\\
\partial_t (\bu_k+\nabla p_{h,k})&\in L^2(0,T;((W_{0}^{1,2}(\hat{\Omega}))^d)^*),
\end{align}
and, for all $\bw\in (W^{1,2}_0(\hat{\Omega}))^d$  and a.a. $t\in (0,T)$, we have
 \begin{equation} \label{eqnew}
     \langle\partial_t (\bu_k+\nabla p_{h,k}),\bw\rangle-\int_{\Omega}\tilde{p}_k\diver \bw+\int_{\Omega} \bD_k:\nabla \bw=-\int_{\Omega}\bH_k:\nabla \bw.
 \end{equation}
Furthermore, these pressures can be constructed so that they satisfy
\begin{align}
  \label{phkstrong2} p_{h,k}&\to 0 && \mbox{strongly in } L^2(0,T;W^{2,2}(\hat{\Omega})),\\
    \label{tildepweak2} \tilde{p}_{k}& \rightharpoonup 0 && \mbox{weakly in } L^2((0,T) \times\hat{\Omega}).
 \end{align}

At this point we are going to use Proposition~\ref{BuLip}. We define
\begin{equation}
\begin{aligned}
\bw_k&:= \bu_k+\nabla p_{h,k},\\
\bG_k&:=\bD_k -\tilde{p}_k\tens{I} + \bH_k
\end{aligned}\label{DFGH}
\end{equation}
and take an arbitrary but fixed set $\tilde{\Omega}$ satisfying the same assumptions as $\Omega$ in \eqref{p000} and $\tilde{\Omega}\subset \overline{\tilde{\Omega}}\subset \hat{\Omega}$. Then, we see that $(\bw_k,\bG_k)$ satisfies the assumptions of Proposition~\ref{BuLip} and, in addition, thanks to \eqref{Hkweak}, \eqref{ukweak2}, \eqref{ukstrong2}, \eqref{phkstrong2} and \eqref{tildepweak2} we see that
\begin{align}
     \bG_k&\rightharpoonup \tens{O} && \mbox{weakly in } (L^2((0,T)\times \tilde{\Omega}))^{d\times d},\label{CC1}\\
     \bw_k&\rightharpoonup \vec{0} && \mbox{weakly in } L^2(0,T; (W^{1,2}(\tilde{\Omega}))^d),\label{CC2}\\
     \bw_k&\to \vec{0}&& \mbox{strongly in } (L^2((0,T)\times\tilde{\Omega}))^d.\label{CC4}
\end{align}
Therefore, for arbitrary $\Lambda>1$, we can find a sequence $\{\bw_{\Lambda,k}\}_{k\in \N}$ for which we can deduce, referring to  \eqref{estnablwlambdak}--\eqref{estwlambdakk} and \eqref{CC2}--\eqref{CC4}, the following convergence results (valid for suitable not relabeled subsequences) as $k\to \infty$:
\begin{align}
     \nabla \bw_{\Lambda,k}&\rightharpoonup^* \tens{O} && \mbox{weakly$^*$ in } (L^{\infty}((0,T)\times \tilde{\Omega}))^{d\times d},\label{CCC1}\\
     \partial_t\bw_{\Lambda,k} \cdot (\bw_{\Lambda,k}-\bw_k) &\rightharpoonup z_{\Lambda} && \mbox{weakly in } L^2((0,T)\times \tilde{\Omega}),\label{CCC2}\\
     \bw_{\Lambda,k}&\to \vec{0} && \mbox{strongly in } (L^2((0,T)\times \tilde{\Omega}))^d.\label{CCC3}
\end{align}
Note that $z_{\Lambda}$ is a general object and we do not know how it is related to the sequence $\bw_k$. However, thanks to the weak lower semicontinuity of the $L^1$ norm and the estimate \eqref{estdtwlambdak}, we can deduce that
\begin{equation}\label{z-est}
\int_0^T \int_{\tilde{\Omega}}|z_{\Lambda}|\le \frac{C}{\sqrt{\Lambda}}.
\end{equation}
Furthermore, we deduce from \eqref{BuLipeq} that, for all $\psi\in \mathcal{C}^{1}_c((0,T)\times\tilde{\Omega})$, there holds
\begin{equation}\label{BuLipeqC}
\begin{split}
    \int_{Q_T} \bG_k: \nabla \bw_{\Lambda,k} \psi= &-\frac{1}{2}\int_{Q_T} (|\bw_{\Lambda,k}|^2-2\bw_k\cdot \bw_{\Lambda,k})\partial_t\psi\\
    &- \int_{Q_T} \bG_k: (\bw_{\Lambda,k}\otimes \nabla \psi)\\
    &-\int_{Q_T}\partial_t \bw_{\Lambda, k}\cdot (\bw_{\Lambda,k}-\bw_k)\psi.
\end{split}
\end{equation}
Recalling the definition of the sets $\mathcal{O}_{\Lambda, k}$, see \eqref{defOLambdak}, and referring to  \eqref{estOLambdak}, we have that 
\begin{align}
\label{estOLambdakC}
|\mathcal{O}_{\Lambda, k}|&\leq \frac{C}{\Lambda}.
\end{align}
In addition, \eqref{estdtwlambdak} and \eqref{CC2} also give the uniform bound
\begin{equation}\label{DD1}
\sup_{\Lambda >1} \, \sup_{k\in \N} \int_0^T \int_{\tilde{\Omega}} |\nabla \bw_{\Lambda,k}|^2+|\nabla \bw_k|^2 <\infty.
\end{equation}

Next, we let $k\to \infty$ in \eqref{BuLipeqC}. Using the convergence results \eqref{CC1}, \eqref{CC4}, \eqref{CCC2} and \eqref{CCC3}, we obtain that, for all $\psi\in \mathcal{C}^{1}_c((0,T)\times\tilde{\Omega})$, there holds
\begin{equation*}
\begin{split}
   \lim_{k\to \infty} \int_{Q_T} \bG_k: \nabla \bw_{\Lambda,k} \psi&= -\int_{Q_T} z_{\Lambda}\psi.
\end{split}
\end{equation*}
Since $\{\bG_k: \nabla \bw_{\Lambda,k}\}_{k\in \N}$ is uniformly bounded in $L^2(Q_T)$, which follows from \eqref{CC1} and \eqref{CCC1}, the above relations says that
\begin{equation}\label{CC10}
\bG_k: \nabla \bw_{\Lambda,k} \rightharpoonup z_{\Lambda} \quad \textrm{ weakly in } L^2((0,T)\times \tilde{\Omega}).
\end{equation}

Next, we show that \eqref{CC10} implies that
\begin{equation}\label{CC11}
\bG_k: \nabla \bw_{k}\rightharpoonup 0 \quad \textrm{ in the biting sense  in } L^1((0,T)\times \tilde{\Omega}).
\end{equation}
By applying Proposition~\ref{biting} and the characterization of $L^1$ weakly convergent sequences to the sequence $\{|\bG_k|^2\}_{k\in \N}$, which is uniformly bounded in $L^1((0,T)\times \tilde{\Omega})$, we get that, for a suitable not relabeled subsequence of $\{|\bG_k|^2\}_{k\in \N}$, there exists a nondecreasing sequence of measurable sets $\{E_j\subset (0,T)\times \tilde{\Omega}\}_{j\in \N}$ fulfilling $|((0,T)\times \tilde{\Omega}) \setminus E_j| \to 0$ as $j\to \infty$ and satisfying that, for any $j\in \N$ and any $\varepsilon >0$, there exists a $\Lambda>1$ such that, for any measurable set $U\subset E_j$ with the property $|U|\le C/\Lambda$, there holds
\begin{equation}
\label{BBC}
\sup_{k\in \N} \int_{U} |\bG_k|^2\le \varepsilon.
\end{equation}
Finally, let $\psi\in L^{\infty}(Q_T)$ be arbitrary. Then, for any $j\in \N$, we have
$$
\begin{aligned}
\lim_{k\to \infty} &\left| \int_{E_j}\bG_k : \nabla \bw_k \psi \right| \\
&\le \lim_{k\to \infty} \left| \int_{E_j}\bG_k : (\nabla \bw_k -\nabla \bw_{\Lambda,k}) \psi \right|+\lim_{k\to \infty} \left| \int_{E_j}\bG_k : \nabla \bw_{\Lambda,k} \psi \right|\\
&=\lim_{k\to \infty} \left| \int_{E_j\cap \mathcal{O}_{\Lambda, k}}\bG_k : (\nabla \bw_k -\nabla \bw_{\Lambda,k}) \psi \right|+ \left| \int_{E_j}z_{\Lambda} \psi \right|\\
&\le C\|\psi\|_{\infty}\left(\frac{1}{\sqrt{\Lambda}} +\limsup_{k\to \infty}  \left(\int_{E_j\cap \mathcal{O}_{\Lambda, k}}|\bG_k|^2\right)^\frac{1}{2}\right),
\end{aligned}
$$
where we have used the convergence result \eqref{CC10}, the definition of the set $\mathcal{O}_{\Lambda, k}$,  the a~priori bounds \eqref{z-est}, \eqref{DD1} and Hölder's inequality. Finally, since the measure of the set $\mathcal{O}_{\Lambda, k}$ tends uniformly to zero as $\Lambda \to \infty$, see \eqref{estOLambdakC}, and also $\{|\bG_k|^2\}_{k\in \N}$ is uniformly equi-integrable on $E_j$, see \eqref{BBC}, we can let $\Lambda \to \infty$ in the above inequality to realize that, for any $\psi \in L^{\infty}(Q_T)$, there holds
$$
\begin{aligned}
\lim_{k\to \infty}\int_{E_j}\bG_k : \nabla \bw_k \psi =0,
\end{aligned}
$$
i.e., we have, for arbitrary $j\in \N$, that
$$
\bG_k : \nabla \bw_k \rightharpoonup 0 \quad \textrm{ weakly in } L^1(E_j).
$$
Since the sequence $\{\bG_k : \nabla \bw_k\}_{k\in \N}$ is uniformly bounded in $L^1((0,T)\times\tilde{\Omega})$, it follows from the properties of~$\{E_j\}_{j\in \N}$ and the above convergence result that we have just recovered~\eqref{CC11}.

Finally, we want to transfer the information obtained in \eqref{CC11} to the original sequence $\{(\tilde{\bG}_k, \tilde{\bw}_k)\}_{k\in \N}$. Using the definitions of $\bG_k$ and $\bw_k$, see \eqref{DFGH} and \eqref{Berlin}, we deduce from \eqref{CC11} that
\small
$$
\begin{aligned}
(\tilde{\bG}_k-\tilde{\bG}-\tilde{p}_k\tens{I} ): \nabla (\tilde{\bw}_k - \tilde{\bw}+\nabla p_{h,k})\rightharpoonup 0 \, \textrm{ in the biting sense in } L^1((0,T)\times \tilde{\Omega}).
\end{aligned}
$$
\normalsize
Thanks to \eqref{bitingC}, \eqref{phkstrong2}, \eqref{tildepweak2}, and the fact that $\tilde{\bw}_k$ and $\tilde{\bw}$ are divergence-free functions, the above relation reduces to
\begin{equation}\label{uztobude}
\begin{aligned}
(\tilde{\bG}_k-\tilde{\bG}): \nabla (\tilde{\bw}_k - \tilde{\bw})\rightharpoonup 0 \quad \textrm{ in the biting sense in } L^1((0,T)\times \tilde{\Omega}).
\end{aligned}
\end{equation}
From  the assumption \eqref{bitingC} it directly follows that
$$
\begin{aligned}
\tilde{\bG}_k: \nabla \tilde{\bw}&\rightharpoonup \tilde{\bG}: \nabla \tilde{\bw} &&\textrm{ weakly in } L^1(Q_T),\\
\tilde{\bG}: \nabla \tilde{\bw}_k&\rightharpoonup \tilde{\bG}: \nabla \tilde{\bw} &&\textrm{ weakly in } L^1(Q_T).
\end{aligned}
$$
Consequently, \eqref{uztobude} directly implies the desired relation
$$
\begin{aligned}
\tilde{\bG}_k: \nabla \tilde{\bw}_k&\rightharpoonup \tilde{\bG}: \nabla \tilde{\bw} &&\textrm{ in the biting sense in } L^1((0,T)\times\tilde{\Omega}).
\end{aligned}
$$
Thus, since $\tilde{\Omega}$ and also $\hat{\Omega}$ with the properties described above were arbitrary and the biting limit of the sequence $\{\tilde{\bG}_k: \nabla \tilde{\bw}_k\}_{k\in \N}$ surely exists on the whole of $Q_T$ (thanks to \eqref{bitingC} the sequence is uniformly bounded in $L^1(Q_T)$), we can use the above identification to conclude \eqref{Br}.

\section{Proof of Theorem~\ref{Burgmain}}\label{Proof2}

This section is devoted to the proof of the main result of the paper. As we have already claimed, we focus only on the three-dimensional setting to simplify the computations. Further, we set the parameters $G,\tau,\rho$ and $2\nu$ to be equal to one, since their size does not play any role in the proof of the existence result as we are dealing with large data analysis. Further, we also set $\vec{f}\equiv \vec{0}$ for simplicity. We would like to emphasize that these simplifications are not at all essential and they do not impact on the validity of the main result.

The proof follows very much the scheme developed in \cite{BuMaLo22}, and therefore, at some parts it is brief. On the other hand, the new parts are proved in full details. In Subsection~\ref{systemwithF}, we introduce a new system of equations involving the matrix $\bF$ and formulate the main theorem concerning the existence of a~solution for such a~system. Note that, as shown in \cite{BuMaLo22} including all details, such a result directly implies the validity of Theorem~\ref{Burgmain}. Then in Subsection~\ref{approxintr}, we introduce an approximative problem, for which the existence theory  directly follows from~\cite{BuMaLo22}. Subsection~\ref{limitk1} is then devoted to the limit procedure based on the weak convergence results. The compactness of the approximations of~$\bF$ is then proved in Subsection~\ref{CFK} and the properties of the determinant of $\bF$ are summarized in Subsection~\ref{limitk2}.

%
%

\subsection{System with evolutionary equation for the tensor \texorpdfstring{$\bF$}{F}}\label{systemwithF}
Referring to \cite{Mas} and \cite{BuMaLo22}, we consider the following evolutionary problem: to find $(\bv,p,\bF):Q_T\to \R^3\times \R \times \R^{3\times 3}$ satisfying 
\begin{align}
    \label{vF1}
    \diver \bv &= 0,\\
    \label{vF2}
    \partial_t \bv + \diver(\bv\otimes \bv) + \nabla p - \diver(\bF \bF^T + \bD) &= \vec{0},\\
    \label{vF3}
    \partial_t \bF + \mbox{Div } (\bF\otimes \bv) - \nabla \bv \bF + \frac{1}{2} (\bF \bF^T \bF - \bF)  &= \tens{O}, \\
    \label{detF}
    \det \bF &>0,
\end{align}
together with the boundary and initial conditions:
\begin{align}
    \label{bcv}
    \bv &= \vec{0} \quad \mbox{on } (0,T)\times\partial\Omega, \\
    \label{icvF} \bv(0, \vec{\cdot}) &= \bv_0 \quad \mbox{and} \quad \bF(0,\vec{\cdot}) = \bF_0 \quad \mbox{in } \Omega.
\end{align}

As one may deduce from \eqref{pepa98}, we have that $\bF\in (L^4(Q_T))^{3\times 3}$ while $\bB=\bF \bF^T\in(L^2(Q_T))^{3\times 3}$. Hence in a weak formulation of \eqref{vF3} the set of admissible test functions is larger than it is for \eqref{vB3lambda}, which is beneficial, for example, for the proof of the compactness of the approximations introduced in the next subsection. The connection between \eqref{vF3} and \eqref{vB3lambda} is achieved as follows: 
if we formally multiply~\eqref{vF3} by~$\bF^T$ from the right, then transpose the equation~\eqref{vF3} and multiply it by~$\bF$ from the left and sum the resulting identities, we arrive at \eqref{vB3lambda}. Although such a procedure is only formal, it may be made rigorous by  repeating step by step the corresponding procedure presented in~\cite{BuMaLo22}. Since such a procedure is independent of dimension, we do not present it here and refer the interested reader to \cite{BuMaLo22}. To conclude, for proving Theorem~\ref{Burgmain} it suffices to prove the following key result.
\begin{thm}
\label{main}
Let $\Omega$ and $T$ satisfy \eqref{p000}. Let  $\bv_0\in L^2_{\bn, \diver}$ and $\bF_0\in (L^2(\Omega))^{3\times 3}$ be such that $\det \bF_0>0$ a.e. in $\Omega$ and $\ln \det \bF_0\in L^1(\Omega)$. Then, there exists a weak solution to \eqref{vF1}--\eqref{icvF}, i.e., there exists a couple $(\bv, \bF)$ satisfying
\begin{align*}
\bv & \in \mathcal{C}_{weak}([0,T]; L^2_{\bn, \diver}) \cap L^2(0,T; W_{\vec{0},\diver}^{1,2}), \\
\partial_t \bv &\in L^{\frac{4}{3}}(0, T; (W_{\vec{0},\diver}^{1,2})^*),\\
\bF & \in \mathcal{C}([0,T]; (L^2(\Omega))^{3\times 3}) \cap (L^4(Q_T))^{3\times3}, \\
\partial_t \bF & \in L^{\frac{4}{3}}(0,T; ((W^{1,2}(\Omega))^{3\times 3})^*),\\
\ln (\det \bF) & \in L^{\infty}(0,T; L^1(\Omega)) \quad \implies \quad \det \bF > 0 \mbox{ a.e. in } Q_T,
\end{align*}
such that, for all $\bw \in W_{\vec{0},\diver}^{1,2}$, $\bA \in (W^{1,2}(\Omega))^{3\times 3}$ and a.a. $t\in (0,T)$, the following identities hold:
\begin{align}
\label{vweak}
\langle \partial_t \bv, \bw \rangle- \int_{\Omega} (\bv \otimes \bv): \nabla \bw + \int_{\Omega} (\bD+\bF \bF^T): \nabla \bw&= 0,\\
\label{Fweak}
\langle \partial_t \bF, \bA \rangle -\int_{\Omega} (\bF \otimes \bv): \nabla \bA - \int_{\Omega} \nabla \bv \bF: \bA + \frac{1}{2}\int_{\Omega} (\bF \bF^T \bF -  \bF): \bA &= 0,
\end{align}
and the initial conditions $\bv_0$ and $\bF_0$ are attained in the following sense:
\begin{equation}
\label{invFweak}
\lim_{t\to 0+} \left(\|\bv(t)-\bv_0\|_2 + \|\bF(t)-\bF_0\|_2\right) = 0.
\end{equation}
\end{thm}

\subsection{Approximation and the first a~priori estimates}
\label{approxintr}
In order to construct the solution according to Theorem~\ref{main}, we consider the following  approximations of the problem \eqref{vF1}--\eqref{icvF}. 

Let $g:\mathbb{R}\to \mathbb{R}$ be a smooth non-negative compactly supported function that satisfies $0\le g(s)\le 1$ for all $s\in \R$,  $g(s)=1$ for all $s\in (-1,1)$ and $g(s)=0$ for all $|s|>2$. Then, for all $k\in \N$, we set $g_k(s):=g(s/k)$ and introduce the following problem: to find $(\bv, p, \bF):Q_T\to \R^3 \times \R \times \R^{3\times 3}$ satisfying
\begin{align}
\label{approxepssigma1}
\diver \bv & = 0,\\
\label{approxepssigma2}
\partial_t \bv + \diver(g_k(|\bv|) \bv\otimes \bv) + \nabla p - \diver(\bD + g_k(|\bF|) \bF\bF^T )& = \vec{0},\\
\label{approxepssigma3}
\partial_t \bF + \mbox{Div } (\bF\otimes \bv) - g_k (|\bF|) (\nabla \bv) \bF + \frac{1}{2} \left(\bF \bF^T \bF - \bF\right) & = \tens{O},\\
\label{detFsigma}\det \bF&>0,
\end{align}
together with the boundary and initial conditions
of the form 
\begin{align}
    \label{bcapproxepssigma}
    &\bv = \vec{0} \quad \mbox{ on } (0,T)\times \partial\Omega, \\
    \label{icapproxepssigma} &\bv(0, \vec{\cdot})  = \bv_0 \quad \mbox{and} \quad \bF(0,\vec{\cdot}) =\bF_0^k:= \bF_0\chi_{\{|\bF_0|\le k\}}+ \bI\chi_{\{|\bF_0|>k\}} \quad \mbox{in } \Omega.   
\end{align}
The problem \eqref{approxepssigma1}--\eqref{icapproxepssigma} can be solved by following the approach and using the tools developed in \cite{BuMaLo22} where the two-dimensional case of the corresponding system without the presence of the factors $g_k(|\bF|)$, $g_k(|\bv|)$ is treated. Indeed here, thanks to the presence of the cut--off function~$g_k$ one has the energy equality at disposal, as in the two-dimensional approach. Hence, we can assume that, for any $k\in \mathbb{N}$, we have a weak solution to \eqref{approxepssigma1}--\eqref{icapproxepssigma}, i.e., there exists $(\bv_k, \bF_k)$ fulfilling
\begin{align*}
\bv_k & \in \mathcal{C}([0,T]; L^2_{\bn, \diver}) \cap L^2(0,T; W_{\vec{0},\diver}^{1,2}), \\
\partial_t \bv_k &\in L^{2}(0, T; (W_{\vec{0},\diver}^{1,2})^*),\\
\bF_k & \in \mathcal{C}([0,T]; (L^2(\Omega))^{3\times 3}) \cap (L^4(Q_T))^{3\times3}, \\
\partial_t \bF_k & \in L^{\frac{4}{3}}(0,T; ((W^{1,2}(\Omega))^{3\times 3})^*),
\end{align*}
such that, for all $\bw \in W_{\vec{0},\diver}^{1,2}$, $\bA \in (W^{1,2}(\Omega))^{3\times 3}$ and a.a. $t\in (0,T)$, the following identities hold ($\bD_k:=\tfrac12(\nabla \bv_k+(\nabla \bv_k)^T)$):
\begin{align}
\label{Mvweak}
&\begin{aligned}
\langle \partial_t \bv_k, \bw \rangle- \int_{\Omega} (g_k(|\bv_k|)\bv_k \otimes \bv_k &) : \nabla \bw \\
&+ \int_{\Omega} (\bD_k+g_k(|\bF_k|)\bF_k \bF_k^T): \nabla \bw= 0,
\end{aligned}\\
\label{MFweak}
&\begin{aligned}
\langle \partial_t \bF_k, \bA \rangle -\int_{\Omega} (\bF_k \otimes \bv_k): \nabla \bA &- \int_{\Omega} g_k(|\bF_k|)\nabla \bv_k \bF_k: \bA \\
&+ \frac{1}{2}\int_{\Omega} (\bF_k \bF^T_k \bF_k -  \bF_k): \bA = 0,
\end{aligned}
\end{align}
and the initial conditions $\bv_0$, $\bF^k_0$ are attained in the following sense
\begin{equation}
\label{MinvFweak}
\lim_{t\to 0+} \left(\|\bv_k(t)-\bv_0\|_2 + \|\bF_k(t)-\bF^k_0\|_2\right) = 0.
\end{equation}

Next, we provide the a~priori estimates for the system \eqref{Mvweak}--\eqref{MFweak}.  This requires us to mollify 
the equation \eqref{MFweak} by means of the Friedrichs lemma on commutators (see \cite{DiPerna1989} or \cite[Theorem 11.2 and Corollary~2.3]{Feireisl2017}). 
Here, we refer for details to  Subsection~\ref{CFK}, Step 1a below, or to \cite{BuMaLo22}, where such a procedure is used many times. 
We set $\bw=\bv_k$ in \eqref{Mvweak} and we take the mollification of $\bF_k$ as a test function in the mollified form of \eqref{MFweak}. Then we sum the resulting identities and, after  integrating   by parts (use that $\bv_k\in L^2(0,T;W_{\vec{0},\diver}^{1,2})$) and over time, using  \eqref{MinvFweak} and taking the limit in the mollifying parameter, arrive, for all $t\in (0,T)$, at
$$
 \|\bv_k(t)\|_2^2 + \|\bF_k(t)\|^2 + \int_0^t\int_{\Omega}2|\bD_k|^2 + |\bF_k \bF_k^T|^2 - |\bF_k|^2=\|\bv_0\|_2^2+\|\bF_0^k\|_2^2.
$$
Taking supremum over $t\in (0,T)$, using Gronwall's lemma, the assumptions on the initial data and Korn's and Young's inequalities, we get
\begin{equation}\label{1sta}
\sup_{t\in (0,T)}\left(\|\bv_k(t)\|_2^2 + \|\bF_k(t)\|_2^2 \right) + \int_0^T \|\nabla \bv_k\|_2^2 + \|\bF_k\|_4^4 \le C(T)\left(\|\bv_0\|_2^2 + |\bF_0\|_2^2 \right).
\end{equation}
Then, it is rather standard to conclude from \eqref{Mvweak}, \eqref{MFweak} and \eqref{1sta},  with the help of Hölder's, Minkowski's and interpolation\footnote{Notice that this is the only  place, which is different when the dimension $d>3$.} inequalities, that
\begin{equation}\label{2sta}
\int_0^T \|\partial_t \bv_k\|_{(W_{\vec{0},\diver}^{1,2})^*}^{\frac{4}{3}} + \|\partial_t \bF_k\|_{((W^{1,2}(\Omega))^{3\times 3})^*}^{\frac43}  + \int_0^T \|\bv_k\|_{\frac{10}{3}}^{\frac{10}{3}} \le C(T,\Omega,\bF_0,\bv_0).
\end{equation}

Finally, 
we need to show that, for all $k\in \N$, $\det \bF_k >0$ a.e. in $Q_T$. This can be done exactly as in \cite[Sect.~6.5]{BuMaLo22} provided that we can justify that $|(\adj \bF_k^T)|$ belongs at least to $L^4(Q_T)$. Since $|(\adj \bF_k^T)|\le C|\bF_k|^{d-1}$, we see that while in two dimensions this requirement follows from the first energy estimate~\eqref{1sta}, in three and higher dimensions, we need to prove  such higher integrability of~$\bF_k$. This is what we do now. In fact, we show that, for any $d\ge 3$ and for all $k\in \N$, we have $\bF_k \in (L^{\infty}(Q_T))^{d\times d}$.  

To prove this, we take as a test function in the mollified form of \eqref{MFweak} the  following function (here, $\bF_k^\delta$ denotes the function $\bF_k$ mollified by the standard kernel with radius $\delta>0$) 
$$
\bA:=\bF_k^\delta \frac{(|\bF_k^\delta|e^{-t/2}-2k)_+}{|\bF_k^\delta|}e^{-t/2}. \qquad \qquad \big(a_+:=\max\{0,a\}\big)
$$
Since $g_k(|\bF_k|)$ vanishes on the set where $|\bF_k|\ge 2k$, we obtain, integrating by parts and over time, taking the limit $\delta\to 0+$ and using \eqref{icapproxepssigma}, the identity 
$$
\begin{aligned}
0&= \lim_{\delta\to0+}\int_0^t\int_{\Omega} \partial_t (|\bF_k^\delta|)(|\bF_k^\delta|e^{-t/2}-2k)_+e^{-t/2} -\frac12|\bF_k^\delta|(|\bF_k^\delta|e^{-t/2}-2k)_+e^{-t/2}\\
&\qquad + \frac{1}{2}\int_0^t\int_{\Omega}|\bF_k\bF_k^T|^2\frac{(|\bF_k|e^{-t/2}-2k)_+}{|\bF_k|}e^{-t/2}\\
&\ge  \lim_{\delta\to0+}\int_0^t\int_{\Omega} \partial_t (|\bF_k^\delta|e^{-t/2})  (|\bF_k^\delta|e^{-t/2}-2k)_{+} =  \frac{1}{2}\int_{\Omega}  (|\bF_k(t)|e^{-t/2}-2k)_+^2
\end{aligned}
$$
valid for all $t\in (0,T)$ (and $k\in \N$ such that $2k\geq \sqrt{d}$).
Hence
\begin{equation}
    |\bF_k|\le 2ke^{T/2} \qquad \textrm{a.e. in } Q_T \qquad \implies \bF_k \in (L^{\infty}(Q_T))^{d\times d}. \label{pepa99}
\end{equation}
Recalling that  $\partial (\ln \det \bA) = \partial \bA: \bA^{-T}$, where $\partial$ stands for $\partial_t$ or $\partial_{x_i}$ with $i=1,...,d$ and $\bF_k^{-T} = (\adj \bF_k^T)/\det \bF_k$, having \eqref{pepa99}, one can now proceed, step by step, as in \cite[Sect.~6.5]{BuMaLo22} and conclude that 
$\det \bF_k >0$ a.e. in~$Q_T$ and the following estimate holds:
\begin{equation}\label{3sta}
\begin{split}
\sup_{t\in (0,T)}\|\ln \det \bF_k(t)\|_1 &\le C(T, \Omega, \|\bF_0^k\|_2, \|\bv_0\|_2,\|\ln \det \bF_0^k\|_1) \\ 
&\le C(T,\Omega,\bF_0,\bv_0).
\end{split}
\end{equation}


\subsection{Weak convergence results as \texorpdfstring{$k\to \infty$}{e0}} \label{limitk1}
The uniform estimates \eqref{1sta} and \eqref{2sta}  imply the existence of $\bv$, $\bF$ such that, for suitably extracted subsequences of $\{\bv_k\}_{k\in \N}$, $\{\bF_k\}_{k\in \N}$ (the subsequences are not relabeled), the following convergence results hold:
\begin{align}
    \label{vkweak}\bv_k&\rightharpoonup \bv && \mbox{weakly-* in } L^\infty(0,T;L^2_{\bn, \diver})\cap L^2(0,T; W_{\vec{0},\diver}^{1,2}),\\
    \label{dtvkweak}\partial_t \bv_k &\rightharpoonup \partial_t \bv && \mbox{weakly in } L^\frac43(0,T; (W_{\vec{0},\diver}^{1,2})^*),\\
    \label{Fkweak}\bF_k&\rightharpoonup \bF && \mbox{weakly-* in } L^\infty(0,T;(L^2(\Omega))^{3\times3})\cap (L^4(Q_T))^{3\times 3},\\
    \label{dtFkweak}\partial_t \bF_k &\rightharpoonup \partial_t \bF&& \mbox{weakly in } L^\frac43(0,T;((W^{1,2}(\Omega))^{3\times3})^*).
\end{align}
It also follows from \eqref{1sta} that 
\begin{equation}
    \label{pepa_a1}
    \begin{split}
    k |\{|\bF_k| > k\}|^{1/4} &\le \left(\int_{\{|\bF_k| > k\}} |\bF_k|^4\right)^\frac14 \le C(T, \bv_0, \bF_0) \\ & \implies \quad |\{|\bF_k|> k\}| \to 0 \textrm{ as } k\to \infty.
    \end{split} 
\end{equation}
Consequently, referring also to the definition of $g_k$, we get (for $k\to \infty$) 
\begin{equation}
    \label{pepa_a2}
    \|1- g_k(|\bF_k|)\|_2^2 = \int_{\{|\bF_k| > k\}} (1-g_k(|\bF_k|))^2 \le |\{|\bF_k| > k\}| \to 0. 
\end{equation}
Applying similar arguments to $\{\bv_k\}_{k\in \N}$ and using basic interpolation inequalities, we conclude that
\begin{align}
\label{gkvkstrong}
    g_k(|\bv_k|)&\to 1 && 
    \mbox{strongly in } L^q(Q_T)\ \forall q \in [1,\infty),\\
    \label{gkFkstrong}
    g_k(|\bF_k|)&\to 1 && 
    \mbox{strongly in } L^q(Q_T)\ \forall q\in [1,\infty).
\end{align}
Next, by \eqref{vkweak}, \eqref{dtvkweak} and the Aubin--Lions compactness lemma we have that
\begin{equation}
    \bv_k\to \bv \quad \mbox{strongly in } (L^q(Q_T))^{3\times3}\ \forall q\in \left[1,\frac{10}{3}\right).\label{vkstrong}
\end{equation}
Moreover, using that $\bv \in L^\infty (0,T;L^2_{\bn,\diver})$, $\partial_t\bv\in L^\frac43 (0,T; (W_{\vec{0},\diver}^{1,2})^*)$, and also $\bF \in L^\infty \left(0,T;(L^2(\Omega))^{3\times3}\right)$ and $\partial_t \bF\in L^\frac{4}{3}(0,T; ((W^{1,2}(\Omega))^{3\times3})^*)$, we see that the functions $\bv$, $\bF$ satisfy
\begin{align}
\label{vcont}
\bv&\in \mathcal{C}_{weak}([0,T]; L^2_{\bn, \diver}), \\
\label{Fcont}
\bF&\in \mathcal{C}_{weak}\left([0,T]; (L^2 (\Omega))^{3\times3}\right).
\end{align}
Finally, as various \emph{nonlinear} (not necessarily tensor)  quantities $\bQ_k$ (such as $\nabla\bv_k \bF_k$, $\bF_k\bF_k^T\bF_k$, $\bF_k\bF_k^T$, $|\bF_k|^2$ and some others introduced below) are uniformly bounded in some Lebesgue spaces $(L^r(Q_T))^{3\times 3}$, $r>1$, there is an element $\overline{\bQ}\in (L^r(Q_T))^{3\times 3}$ and a subsequence of $\bQ_k$ (denoted again $\bQ_k$) so that 
\begin{equation}
\label{conven}
\bQ_k \rightharpoonup \overline{\bQ} \textrm{ weakly in } (L^r(Q_T))^{3\times 3} \quad \textrm{ as } k\to \infty. 
\end{equation}
We do not list all such quantities here, but we use this notation and convention (that all these selections of weakly converging sequences are made at this moment) in what follows.
Referring to this convention and applying the above convergence results \eqref{vkweak}--\eqref{vkstrong} to the identities~\eqref{Mvweak} and~\eqref{MFweak}, we obtain, for all $\bw\in W_{\vec{0},\diver}^{1,2}$, $\bA\in (W^{1,2}(\Omega))^{3\times 3}$ and a.a. $t\in (0,T)$, that
\begin{align} \label{2where}\langle\partial_t \bv, \bw\rangle - \int_{\Omega} (\bv \otimes \bv): \nabla\bw + \int_{\Omega} \bD: \nabla \bw + \int_{\Omega} \overline{\bF \bF^T}: \nabla \bw& = 0, \\
\label{Fwithoutcomp} \langle\partial_t \bF, \bA\rangle - \int_{\Omega}(\bF \otimes\bv): \nabla\bA - \int_{\Omega} \overline{(\nabla \bv) \bF}:\bA + \frac{1}{2}\int_{\Omega}(\overline{\bF \bF^T \bF} -\bF):\bA &= 0.
\end{align}

In order to show that $\bF\in \mathcal{C}([0,T];(L^2(\Omega))^{3\times3})$ and \eqref{invFweak} hold, one can proceed step by step as in the proofs of the corresponding properties in the two-dimensional setting, see \cite[Sects. 6.3 and B.4]{BuMaLo22}. Hence, we do not present the proofs of these results here. We focus on the identification of the weak limits in the nonlinear terms, i.e.  on showing that $\overline{\bF \bF^T} = \bF \bF^T$, $\overline{\nabla \bv \bF} = \nabla \bv \bF$ and $\overline{\bF \bF^T \bF}= \bF \bF^T \bF$ in \eqref{2where} and \eqref{Fwithoutcomp}. Thanks to \eqref{vkweak} and \eqref{Fkweak} it is sufficient, for the above identifications, to prove the compactness of $\{\bF_k\}_{k\in \N}$ in $(L^2 (Q_T))^{3\times3}$.  Last but not least, we prove that $\det \bF>0$ a.e. in $Q_T$, which is, however, a simple consequence of the compactness of $\{\bF_k\}_{k\in \N}$ in $(L^2(Q_T))^{3\times3}$ and the uniform bound~\eqref{3sta}.

\subsection{Compactness of \texorpdfstring{$\{\bF_k\}_{k\in \N}$}{Fs} in \texorpdfstring{$(L^2(Q_T))^{3\times3}$}{LQ}}\label{CFK}
In principle we follow the proof of the compactness of $\{\bF_\varepsilon\}$ in $(L^2(Q_T))^{2\times2}$ performed in~\cite[Sect.~6.4]{BuMaLo22}. This means that we work with the difference between the equation for $\bF_k$ formally multiplied by $\bF_k$ and the equation for $\bF$ formally multiplied by $\bF$. After computations, in which we use, among other things, the evolutionary equation for $\bv_k$ and the convergence results associated with the sequence $\{\bv_k\}_{k\in \N}$, we arrive at the inequality rigorously written, for all non-negative $\varphi\in \mathcal{C}_c^\infty((-\infty,T)\times\Omega)$,  as
\begin{equation}\label{prepGron}
    -\int_{Q_T}(\overline{|\bF|^2}-|\bF|^2)\partial_t \varphi-\int_{Q_T} (\overline{|\bF|^2}-|\bF|^2) \bv\cdot \nabla \varphi\leq \int_{Q_T} L(\overline{|\bF|^2}-|\bF|^2)\varphi,
\end{equation}
where $L$ is an $L^2(Q_T)$ function specified below. Due to the weak lower semicontinuity of $L^2$ norm, we also have that $\overline{|\bF|^2}\ge |\bF|^2$ a.e. in $Q_T$.  Then, we can employ Proposition~\ref{convae} with $f:=\overline{|\bF|^2}-|\bF|^2$ and  conclude that
\begin{equation}\label{compae}
    \overline{|\bF|^2}=|\bF|^2 \quad \mbox{a.e. in } Q_T,
\end{equation}
which is equivalent to the compactness of $\{\bF_k\}_{k\in \N}$ in $(L^2(Q_T))^{3\times3}$, see \cite{BuMaLo22}.

Thus, in what follows, we focus on the derivation of \eqref{prepGron} and this is indeed the key difference compared with the procedure  presented in \cite{BuMaLo22} as $\bv$ is not an admissible test function in \eqref{2where} (which, after integration over $(0,T)$, holds for every $\bw\in L^4(0,T; W_{\vec{0},\diver}^{1,2})$, while $\bv$ belongs only to $L^2(0,T; W_{\vec{0},\diver}^{1,2})$) and consequently one does not have the energy equality for the velocity at disposal. We split the proof into the following steps.

\bigskip
\noindent \textbf{Step 1a:} In the first step we derive the equations for $|\bF|^2$ and also the renormalized equation for $|\bF_k|^2$. This means that we want to show that, for any $m,k\in \N$ and  for all $\varphi\in C_c^\infty((-\infty,T)\times\Omega)$, there holds
\begin{equation}
    \label{step1Fsigmaa}
    \begin{split}
       & -\int_{Q_T} \mathcal{G}_m(|\bF_k|^2)\partial_t \varphi-\int_{\Omega} \mathcal{G}_m(|\bF_0^k|^2)\varphi(0)-\int_{Q_T} \mathcal{G}_m(|\bF_k|^2)\bv_k\cdot \nabla \varphi\\
        &\qquad -2\int_{Q_T} (\nabla \bv_k\bF_k): g_m(|\bF_k|^2)g_k(|\bF_k|)\bF_k\varphi\\
        &\qquad +\int_{Q_T}g_m(|\bF_k|^2)\left(|\bF_k \bF_k^T|^2-|\bF_k|^2\right)\varphi
        = 0
    \end{split}
\end{equation}
and
\begin{equation}
\label{step1Fsigmab}
\begin{split}
    -\int_{Q_T} |\bF|^2\partial_t \varphi-\int_{\Omega} |\bF_0|^2\varphi(0)-\int_{Q_T}|\bF|^2 \bv\cdot \nabla\varphi-2\int_{Q_T}\overline{\nabla \bv \bF}: \bF\varphi\\+\int_{Q_T}\left(\overline{\bF \bF^T \bF}:\bF-|\bF|^2\right)\varphi=0.
    \end{split}
\end{equation}
Here, recalling that $g$ is introduced at the beginning of Subsect.~\ref{approxintr},  the function $\mathcal{G}_m$ is defined as
\begin{equation}
\label{defGms}
\mathcal{G}_m(s):=\int_0^s g_m(u)\, \mathrm{d}u, \quad\textrm{ where } g_m(s)= g\left(\frac{s}{m}\right).
\end{equation}

\bigskip
\noindent \textbf{Step 1b:}
We subtract \eqref{step1Fsigmab} from \eqref{step1Fsigmaa},  take the limit $k\to\infty$, and then the limit $m\to \infty$, and we get the following inequality, valid for all non-negative  $\varphi\in \mathcal{C}_c^\infty((-\infty,T)\times\Omega)$:
\begin{equation}
    \label{step1compFsigma}
    \begin{split}
    -&\int_{Q_T} \left(\overline{|\bF|^2}-|\bF|^2\right)\partial_t \varphi - \int_{Q_T} \left(\overline{|\bF|^2}-|\bF|^2\right) \bv\cdot \nabla \varphi \\
    &\qquad \leq \int_{Q_T} \left(\overline{|\bF|^2}-|\bF|^2\right)\varphi-\ell_{1}+2\ell_{2},
\end{split}
\end{equation}
where
\begin{align}
\label{lb}
 \ell_{1}&:= \lim_{m\to \infty}\lim_{k\to \infty} \int_{Q_T} \left(|\bF_k \bF_k^T|^2g_m(|\bF_k|^2)-\bF_k\bF_k^T\bF_k:\bF\right)\varphi,
    \\
\label{la}
    \ell_{2}&:= \lim_{m\to \infty}\lim_{k\to \infty} \int_{Q_T} \left(\nabla \bv_k: g_m(|\bF_k|^2) g_k(|\bF_k|)\bF_k \bF_k^T- \nabla\bv_k\bF_k:\bF\right)\varphi.
\end{align}

\bigskip

\noindent \textbf{Step 2a:} We prove that $\ell_{1}\ge 0$ for all non-negative $\varphi\in \mathcal{C}_c^\infty((-\infty,T)\times\Omega)$.

\bigskip

\noindent \textbf{Step 2b:} We prove that, for all non-negative $\varphi\in \mathcal{C}_c^\infty((-\infty,T)\times\Omega)$, there holds
\begin{equation}
\label{step2compFsigma}
    \ell_{2}\leq \int_{Q_T} \tilde{L} \left(\overline{|\bF|^2}-|\bF|^2\right)\varphi,
\end{equation}
where $\tilde{L}$ is an $L^2(Q_T)$ function specified below.


\bigskip

\noindent \textbf{Step 3:} Here we summarize the results from Step~1 and Step~2 and deduce the inequality \eqref{prepGron}, which suffices to the proof.

\bigskip

\begin{proof}[\bf Ad Step 1a:]  Let us extend $\bv_k,\bv$ by zero and $\bF_k,\bF$ by the identity tensor outside of $\Omega$ and use $h^\delta$ to denote the mollification of $h\in L^1_{loc}(\R^3)$ by the standard mollifying kernel $\omega_\delta$ ($\delta>0$) with respect to the spatial variable.

We prove \eqref{step1Fsigmaa}  first.
In \eqref{MFweak}  we set \[\bA=\bA(\vec{\cdot}):= \omega_\delta(\vec{y}-\vec{\cdot})2g_m(|\bF_k^\delta(\vec{y})|^2)\bF_k^\delta(\vec{y}),\] where $\vec{y}\in\Omega$ and $m\in\N$ are arbitrary.
As $\bF_k\in (L^4(Q_T))^{3\times3}$ and $\bv_k\in L^2(0,T; W_{\vec{0},\diver}^{1,2})$, we obtain (for more details see \cite{BuMaLo22}) that
$
    \partial_t \bF_k^\delta\in L^\frac{4}{3} \left(0,T; \mathcal{C}^1 (\overline{\Omega})\right)
$
and, as also $\diver\bv_k=0$, we have a.e. in $Q_T$, that
\begin{equation}
\label{Fsigmaalphaae}
\begin{split}
    \partial_t \mathcal{G}_m(|\bF_k^\delta|^2)+\diver (\mathcal{G}_m(|\bF_k^\delta|^2)\bv_k) - 2(g_k(|\bF_k|)\nabla \bv_k\bF_k)^\delta: g_m(|\bF_k ^\delta|^2)\bF_k^\delta\\
    +\left((\bF_k \bF_k^T \bF_k)^\delta:\bF_k^\delta-|\bF_k^\delta|^2\right)g_m(|\bF_k^\delta|^2) 
    =2\tens{E}_k^\delta: g_m(|\bF_k^\delta|^2)\bF_k^\delta,
    \end{split}
\end{equation}
where
$
    \tens{E}_k^\delta:= \mbox{Div } (\bF_k^\delta\otimes \bv_k)-\mbox{Div } (\bF_k\otimes\bv_k)^\delta.
$
We multiply \eqref{Fsigmaalphaae} by $\varphi\in \mathcal{C}_c^\infty((-\infty,T)\times\Omega)$, integrate the result over $Q_T$ and let $\delta\to 0_+$. Employing integration by parts, the facts that $\bF_k^\delta, \bF_k\in \mathcal{C}([0,T]; (L^2(\Omega))^{3\times3})$ and $\bF_k(0)=\bF_0^k$ a.e. in $\Omega$, the property $\bv_k\in L^2(0,T;W_{\vec{0},\diver}^{1,2})$, the Friedrichs commutator lemma (applied on $\tens{E}_k^\delta$ with $\bv_k~\in ~L^2(0,T;W_{\vec{0},\diver}^{1,2})$ and $\bF_k\in (L^4(Q_T))^{3\times3}$) together with Hölder's inequality 
and standard properties of mollifying, we obtain \eqref{step1Fsigmaa}.

Using very similar arguments, we observe, setting $\bA(\vec{\cdot}):= \omega_\delta(\vec{y}-\vec{\cdot}) \bF^\delta (\vec{y})\varphi(\vec{y})$ in \eqref{Fwithoutcomp}, integrating the result with respect to $[t,\vec{y}]$ over $Q_T$ and letting $\delta\to0+$, that \eqref{step1Fsigmab} holds. Note that this procedure is described in detail in ~\cite{BuMaLo22}.
\end{proof}

\begin{proof}[\bf Ad Step 1b]
Let us subtract \eqref{step1Fsigmab} from \eqref{step1Fsigmaa}, take the limit $k\to\infty$, and then let $m\to \infty$. Note that, by the definitions of $\bF_0^k$ and $\mathcal{G}_m$ (see \eqref{icapproxepssigma} and \eqref{defGms}), $\mathcal{G}_m(|\bF_0^k|^2)\to \mathcal{G}_m(|\bF_0|^2)$ strongly in $L^1(\Omega)$ as $k\to \infty$ and $\mathcal{G}_m(|\bF_0|^2)\to |\bF_0|^2$ strongly in $L^1(\Omega)$ as $m\to\infty$. Hence, 
to conclude \eqref{step1compFsigma} it  suffices to show, for all non-negative $\varphi\in \mathcal{C}_c^\infty((-\infty,T)\times\Omega)$, the following equalities
\begin{equation}
    \label{limsigmaalpha1}
    \begin{aligned}
    &\lim_{m\to \infty} \lim_{k\to\infty} \int_{Q_T}\left(\mathcal{G}_m(|\bF_k|^2)\partial_t \varphi+ \mathcal{G}_m(|\bF_k|^2)\bv_k\cdot \nabla \varphi\right)
    \\
    &\qquad= \int_{Q_T} \left(\overline{|\bF|^2}\partial_t \varphi +
    \overline{|\bF|^2}\bv\cdot \nabla \varphi\right)
    \end{aligned}
\end{equation}
and
\begin{equation}
    \label{limsigmaalpha2}
    \lim_{m\to \infty} \lim_{k\to\infty} \int_{Q_T} g_m(|\bF_k|^2)|\bF_k|^2\varphi= \int_{Q_T} \overline{|\bF|^2} \varphi.
\end{equation}

To prove the above identities, we first observe
 that
\begin{equation}\label{estgmFk}
\begin{aligned}
&\mathcal{G}_m(|\bF_k|^2) \le |\bF_k|^2, \qquad |\bF_k|^2 -\mathcal{G}_m(|\bF_k|^2)\le |\bF_k|^2 \chi_{\{|\bF_k|^2>m\}},\\
&g_m(|\bF_k|^2)|\bF_k|^2\le |\bF_k|^2, \qquad    |\bF_k|^2-g_m(|\bF_k|^2)|\bF_k|^2 \le |\bF_k|^2 \chi_{\{|\bF_k|^2>m\}}.
\end{aligned}
\end{equation}
Therefore, we can use the uniform boundedness of $\{\bF_k\}_{k\in \N}$ in $(L^4(Q_T))^{3\times3}$ and the convention \eqref{conven} 
and conclude that, as $k \to \infty$, we have the following
\begin{equation}\label{S1a0}
\begin{aligned}
\mathcal{G}_m(|\bF_k|^2)& \rightharpoonup \overline{\mathcal{G}_m(|\bF|^2)} &&\textrm{weakly in }L^2(Q_T),\\
g_m(|\bF_k|^2)|\bF_k|^2 & \rightharpoonup \overline{g_m(|\bF|^2)|\bF|^2} &&\textrm{weakly in }L^2(Q_T).
\end{aligned}
\end{equation}
In addition, using the weak lower semicontinuity of the~$L^2$-norm, we have that
$$
\begin{aligned}
\int_{Q_T} \overline{\mathcal{G}_m(|\bF|^2)}^2 + \overline{g_m(|\bF|^2)|\bF|^2}^2 &\le \liminf_{k\to \infty}\int_{Q_T} (\mathcal{G}_m(|\bF_k|^2))^2 + (g_m(|\bF_k|^2)|\bF_k|^2)^2\\
&\le 2\liminf_{k\to \infty}\int_{Q_T}|\bF_k|^4 \le C,
\end{aligned}
$$
with $C$ being independent of $m$. Consequently, there exist functions $\overline{\mathcal{G}_{\infty}(|\bF|^2)}$ and $\overline{g_{\infty}(|\bF|^2)|\bF|^2}$ fulfilling, as $m\to \infty$, that (modulo subsequences)
\begin{equation}\label{S1a1}
\begin{aligned}
\overline{\mathcal{G}_m(|\bF|^2)}& \rightharpoonup \overline{\mathcal{G}_{\infty}(|\bF|^2)} &&\textrm{weakly in }L^2(Q_T),\\
\overline{g_m(|\bF|^2)|\bF|^2} & \rightharpoonup \overline{g_{\infty}(|\bF|^2)|\bF|^2} &&\textrm{weakly in }L^2(Q_T).
\end{aligned}
\end{equation}
To identify the limiting objects $\overline{\mathcal{G}_{\infty}(|\bF|^2)}$ and $\overline{g_{\infty}(|\bF|^2)|\bF|^2}$, we can use the weak lower semicontinuity  of the~$L^1$-norm and the estimates \eqref{estgmFk} and \eqref{1sta} to obtain that
$$
\begin{aligned}
&\int_{Q_T}\left|\overline{|\bF|^2}-\overline{\mathcal{G}_m(|\bF|^2)}\right| + \left|\overline{|\bF|^2}-\overline{g_m(|\bF|^2)|\bF|^2}\right| \\
&\quad \le \liminf_{k\to \infty}  \int_{Q_T}\left||\bF_k|^2-\mathcal{G}_m(|\bF_k|^2)\right| + \left||\bF_k|^2-g_m(|\bF_k|^2)|\bF_k|^2\right|\\
&\quad \le 2\liminf_{k\to \infty}  \int_{Q_T\cap\{|\bF_k|^2>m\}} |\bF_k|^2\le 2\liminf_{k\to \infty} \int_{Q_T}\frac{|\bF_k|^4}{m} \le \frac{C}{m}.
\end{aligned}
$$
Thus, letting $m\to \infty$, we see that
\begin{equation}\label{S1a2}
\begin{aligned}
\overline{\mathcal{G}_m(|\bF|^2)}& \to \overline{|\bF|^2} &&\textrm{strongly in }L^1(Q_T),\\
\overline{g_m(|\bF|^2)|\bF|^2} & \to \overline{|\bF|^2} &&\textrm{strongly in }L^1(Q_T).
\end{aligned}
\end{equation}
Thanks to the uniqueness of the weak limits in \eqref{S1a1} we have
\begin{equation}\label{S1a3}
\overline{\mathcal{G}_{\infty}(|\bF|^2)}=\overline{g_{\infty}(|\bF|^2)|\bF|^2}=\overline{|\bF|^2}.
\end{equation}
Hence, to justify \eqref{limsigmaalpha1} and \eqref{limsigmaalpha2},  we use \eqref{S1a0} and \eqref{vkstrong} to identify the limit with respect to $k\to \infty$, then use \eqref{S1a1} and the fact that $\bv\in (L^2(Q_T))^3$ to identify the limit as $m\to \infty$ and finally use \eqref{S1a3}. This completes Step~1b.
\end{proof}

\begin{proof}[\bf Ad Step 2a] To show the non-negativity of $\ell_{1}$, we can rewrite this term and we see that what we need to prove is the following inequality
\begin{equation}
\label{monsigmaalphapomMM}
\begin{split}
    \lim_{m\to \infty} \lim_{k\to\infty} &\left(\int_{Q_T} \left(|\bF_k \bF_k^T|^2-\bF_k \bF_k^T \bF_k:\bF\right)g_m(|\bF_k|^2)\varphi\right.\\
    &\qquad \left.+\int_{Q_T}\left(g_m(|\bF_k|^2)\bF_k \bF_k^T \bF_k - \bF_k \bF_k^T \bF_k\right):\bF\varphi \right)\geq 0.
\end{split}
\end{equation}
The first term of \eqref{monsigmaalphapomMM} can be written as
\begin{equation}
\label{monsigmaalphapom2}
\begin{split}
    \lim_{m\to \infty} \lim_{k\to\infty} &\int_{Q_T}\left((\bF_k \bF_k^T \bF_k -\bF \bF^T \bF):(\bF_k-\bF)g_m(|\bF_k|^2)\right.\\
    &\qquad \left.+\bF \bF^T \bF: (\bF_k-\bF)g_m(|\bF_k|^2)\right)\varphi.
\end{split}
\end{equation}
Using the fact that the matrix function $S(\bA):=\bA \bA^T \bA$ is monotone, see \cite[Lemma 4.2]{BuMaLo22}, and also the non-negativity of $g_m$ and $\varphi$, we see that the first term in \eqref{monsigmaalphapom2} is non-negative. Hence, to prove~\eqref{monsigmaalphapomMM} it suffices to show that
\begin{equation*}
\begin{split}
    \lim_{m\to \infty} \lim_{k\to\infty} &\left(\int_{Q_T} \left(g_m(|\bF_k|^2)\bF_k \bF_k^T \bF_k - \bF_k \bF_k^T \bF_k\right):\bF\varphi\right.\\
    &\qquad +\left.\int_{Q_T}\bF \bF^T \bF:(\bF_k-\bF)g_m(|\bF_k|^2)\varphi \right)=0.
\end{split}
\end{equation*}
However, the weak convergence result \eqref{Fkweak} implies that the above relation is equivalent to 
\begin{equation}
\label{monsigmaalphapom3}
    \lim_{m\to \infty} \lim_{k\to\infty} \int_{Q_T} \left(\bF_k \bF_k^T \bF_k :\bF+\bF \bF^T \bF:(\bF_k-\bF)\right)\varphi (g_m(|\bF_k|^2)-1)=0.
\end{equation}
Using the definition of $g_m$ and the estimate \eqref{1sta}, we deduce that
$$
\begin{aligned}
&\left|\int_{Q_T} \left(\bF_k \bF_k^T \bF_k :\bF+\bF \bF^T \bF:(\bF_k-\bF)\right)\varphi (g_m(|\bF_k|^2)-1)\right|\\
&\qquad \le \|\varphi\|_{L^{\infty}(Q_T)} \int_{\{|\bF_k|^2>m\}} |\bF_k|^3 |\bF| + |\bF|^4+ |\bF|^3 |\bF_k|\\
&\qquad \le C \left(\int_{\{|\bF_k|^2>m\}} |\bF|^4 +\left(\int_{\{|\bF_k|^2>m\}} |\bF|^4\right)^{\frac34}+\left(\int_{\{|\bF_k|^2>m\}} |\bF|^4\right)^{\frac14} \right).
\end{aligned}
$$
Since 
\begin{equation}\label{measFm}\left|\{|\bF_k|^2>m\}\right|\le \int_{\{|\bF_k|^2> m\}}\frac{|\bF_k|^4}{m^2}\le \frac{C}{m^2} \to 0 \textrm{ as } m\to \infty,
\end{equation} 
and $\bF\in (L^4(Q_T))^{3\times 3}$, we observe that
$$
\sup_{k\in \N}\int_{\{|\bF_k|^2>m\}} |\bF|^4 \to 0 \textrm{ as } m\to \infty.
$$
Consequently, the relation~\eqref{monsigmaalphapom3} holds, which completes Step 2a.
\end{proof}


\begin{proof}[\bf Ad  Step 2b]
We use the symmetry of $\bF \bF^T$ and $\bF_k \bF^T_k$ and  decompose $\ell_{2}$ (see \eqref{la}) into the sum $\sum_{\gamma=1}^3 L_\gamma$,
where
\begin{align}
    L_1&:= \lim_{m\to \infty} \lim_{k\to \infty}\int_{Q_T} \left(\bD_k: g_m(|\bF_k|^2)g_k(|\bF_k|)\bF_k \bF_k^T-\bD:g_k(|\bF_k|)\bF_k \bF_k^T\right)\varphi,\label{L1}\\
    L_2&:= \lim_{k\to\infty}\int_{Q_T} \left(\bD:(g_k(|\bF_k|)\bF_k \bF_k^T - \bF \bF^T)\right)\varphi,\label{L2}\\
    L_3&:= \lim_{k\to\infty}\int_{Q_T} \left(\nabla \bv\bF:\bF- \nabla \bv_k \bF_k: \bF\right)\varphi.\label{L3}
\end{align}
Recall that
$$
\bD_k:= \frac12 \left( \nabla \bv_k + (\nabla \bv_k)^T \right), \qquad \bD:= \frac12 \left( \nabla \bv + (\nabla \bv)^T \right).
$$

We start with the easiest term $L_2$.  Using \eqref{Fkweak}, \eqref{gkFkstrong} and the facts that $\bD$ is a fixed function from $(L^2(Q_T))^{3\times3}$ and that $\varphi$ is non-negative and bounded, we observe that
\begin{equation*}
    L_2 =\lim_{k\to\infty} \int_{Q_T} \bD :(g_k(|\bF_k|)(\bF_k-\bF)(\bF_k^T-\bF^T))\varphi,
\end{equation*}
and consequently 
\begin{equation}
    \label{step2Fsigmaprep}
    \begin{split}
    |L_2|&\leq \lim_{k\to\infty} \int_{Q_T} |\bD| |\bF_k-\bF|^2 \varphi \\ &=\lim_{k\to \infty} \int_{Q_T} |\bD|(|\bF_k|^2-|\bF|^2)\varphi= \int_{Q_T} |\bD| \left(\overline{|\bF|^2}-|\bF|^2\right) \varphi.
    \end{split}
\end{equation}

In what follows, we shall treat the terms $L_1$ and $L_3$. To be able to perform the required computations, we first decompose the velocity fields as 
    $$\bv_k:=\bv_k^1+\bv_k^2,$$
where $\bv_{k}^1$ is the weak solution $\bu_k$ to the Stokes problem \eqref{dva12}--\eqref{dva13} with $\bu_{0,k}:=\bv_0$ (for all $k\in \N$) and with
\begin{equation}
\label{defHkconcl1}
\bH_k:=g_k(|\bF_k|)\bF_k \bF^T_k.
\end{equation}
Note that, by \eqref{1sta}, the sequence $\{\bH_k\}_{k\in \N}$ defined by \eqref{defHkconcl1} is uniformly bounded in $(L^2(Q_T))^{3\times 3}$, and that, thanks to Proposition~\ref{recpres1}, the velocities $\bv_{k}^1$ are uniquely defined. In addition, thanks to \eqref{ukweak1} and \eqref{dtukweak},  there exists $\bv^1$ such that
\begin{align}
\label{vk1weak}
\bv_{k}^1&\rightharpoonup \bv^1 &&\textrm{ weakly in }L^2(0,T; W^{1,2}_{0,\diver}),\\
\label{dtvk1weak}
\partial_t \bv_{k}^1&\rightharpoonup \partial_t \bv^1 &&\textrm{ weakly in }L^2(0,T; (W^{1,2}_{0,\diver})^*).
\end{align}
Since $\bv_k^2= \bv_k-\bv_k^1$, we see that $\bv_k^2$ is the solution $\bu_k$ to the equation \eqref{SII} with the initial condition $\bu_k(0):=\bv_k^2(0)=\vec{0}$ and with 
\begin{equation}
\label{defHkconcl2}
\bH_k:=- g_k(|\bv_k|)(\bv_k\otimes \bv_k).
\end{equation}
By \eqref{vkstrong} and \eqref{gkvkstrong}, the sequence $\{\bH_k\}_{k\in \N}$ defined by \eqref{defHkconcl2} converges to $\bH:=-(\bv\otimes \bv)$ strongly in $(L^q (Q_T))^{3\times 3}$ for all $q\in [1,5/3)$,
and, obviously, thanks to \eqref{vkweak} and \eqref{vk1weak}, we have that
\begin{equation}
\label{vk2weak}
\bv_{k}^2\rightharpoonup \bv^2 \textrm{ weakly in }L^2(0,T; W^{1,2}_{0,\diver}),
\end{equation}
where $\bv^2=\bv-\bv^1$.
Proposition \ref{recpres2} then implies that 
\begin{equation}
\label{vk2strong}
\nabla\bv_{k}^2\to \nabla \bv^2 \quad \textrm{ strongly in }(L^r(Q_T))^{3\times3} \textrm{ for all } r\in [1,2).
\end{equation}
Next, let us denote
$$
\bD_k^i:=\frac12 \left( \nabla \bv_{k}^i +(\nabla \bv_{k}^i)^T\right), \quad \bD^i:=\frac12 \left( \nabla \bv^{i} +(\nabla \bv^{i})^T\right), \quad i=1,2.
$$
We shall deal with biting limits of terms involving $\bD^1_k$, hence we also define the 
sequence $\{b_k\}_{k\in \N}$ as
\begin{equation}\label{CISLO}
b_k:= |\bD_k^1|^2+|\bF_k|^4. 
\end{equation}
Thanks to \eqref{1sta} and \eqref{vk1weak}, the sequence $\{b_k\}_{k\in \N}$ is uniformly bounded in $L^1(Q_T)$. 
Hence, Chacon's biting lemma (Proposition~\ref{biting}) and the characterisation of $L^1$-weakly converging sequences imply, for a suitable not relabeled subsequence of $\{b_k\}_{k\in \N}$, the existence of a nondecreasing sequence of measurable sets $\{E_j\}_{j\in \N}$ ($E_j\subset Q_T$ for all $j\in \N$) with the property $|Q_T \setminus E_j|\to 0$ as $j\to \infty$ and satisfying that, for any $j\in \N$ and any $\varepsilon>0$, there exists $\delta>0$ such that, for any measurable set $U\subset E_j$ with $|U|<\delta$, one has
\begin{equation}
\sup_{k\in \N} \int_{U} b_k \le \varepsilon.\label{EQI}
\end{equation}
Now we decompose the terms $L_1$ and $L_3$ in the following way:
\begin{equation}
\label{L13dec}
L_1=L_1^{j,1}+L_1^{\tilde{j},1}+L_1^{2}, \qquad L_3=L_3^{j,1}+L_3^{\tilde{j},1}+L_3^{2},
\end{equation}
where 
\begin{align*}
    &L_1^{j,1}:= \lim_{m\to \infty} \lim_{k\to \infty}\int_{E_j} \left(\bD^1_k: g_m(|\bF_k|^2)g_k(|\bF_k|)\bF_k \bF_k^T-\bD^1:g_k(|\bF_k|)\bF_k \bF_k^T\right)\varphi,\\
     &L_1^{\tilde{j},1}:= \lim_{m\to \infty} \lim_{k\to \infty}\int_{Q_T\setminus E_j} \left(\bD^1_k: g_m(|\bF_k|^2)g_k(|\bF_k|)\bF_k \bF_k^T-\bD^1:g_k(|\bF_k|)\bF_k \bF_k^T\right)\varphi,\\
     &L_1^2:=\lim_{m\to \infty} \lim_{k\to \infty}\int_{Q_T} \left(\bD^2_k: g_m(|\bF_k|^2)g_k(|\bF_k|)\bF_k \bF_k^T-\bD^2:g_k(|\bF_k|)\bF_k \bF_k^T\right)\varphi,\\
    &L_3^{j,1}:= \lim_{k\to\infty}\int_{E_j} \left(\nabla \bv^1\bF:\bF- \nabla \bv^1_k \bF_k: \bF\right)\varphi,\\
    &L_3^{\tilde{j},1}:= \lim_{k\to\infty}\int_{Q_T\setminus E_j} \left(\nabla \bv^1\bF:\bF- \nabla \bv^1_k \bF_k: \bF\right)\varphi,\\
    &L_3^{2}:= \lim_{k\to\infty}\int_{Q_T} \left(\nabla \bv^2\bF:\bF - \nabla \bv^2_k \bF_k: \bF\right)\varphi.
\end{align*}

First, we focus on the limits $L_1^{j,1}$ and $L_3^{j,1}$. We start with an identification of $L_1^{j,1}$. Using the definition of $g_m$ (resp. $g_k$), Young's inequality, the estimate \eqref{measFm} and the uniform equi-integrability of the sequence $\{b_k\}_{k\in \N}$ on $E_j$ (see \eqref{EQI} and recall \eqref{CISLO} for the definition of $b_k$), we observe that
$$
\begin{aligned}
&\lim_{m\to \infty} \lim_{k\to \infty}\left|\int_{E_j} \left(\bD^1_k: (g_m(|\bF_k|^2)-1)g_k(|\bF_k|)\bF_k \bF_k^T\right)\varphi\right|\\
&\quad \le \|\varphi\|_{L^\infty(Q_T)}\lim_{m\to \infty} \lim_{k\to \infty}\int_{E_j \cap \{|\bF_k|^2>m\}} b_k=0.
\end{aligned}
$$
 Therefore, we may simplify the definition of $L_1^{j,1}$ to
$$
L_1^{j,1}= \lim_{k\to \infty}\int_{E_j} \left(\bD^1_k: g_k(|\bF_k|)\bF_k \bF_k^T-\bD^1:g_k(|\bF_k|)\bF_k \bF_k^T\right)\varphi.
$$
To identify the last limit, we employ Theorem~\ref{bitingT} with $\tilde{\bw}_k:= \bv^1_k$, $\tilde{\bw}:= \bv^1$, $\tilde{\bG}_k:= g_k(|\bF_k|) \bF_k \bF_k^T+\bD^1_k$ and $\tilde{\bG}:= \overline{\bF \bF^T}+\bD^1$. Note that the estimate \eqref{EQI} implies the uniform equi-integrability of the sequence $\{\tilde{\bG}_k:\bD_k^1\}_{k\in \N}$ on $E_j$. Due to the convergences \eqref{vk1weak}, \eqref{dtvk1weak} 
and \eqref{gkFkstrong} (together with the uniform boundedness of $\{\bF_k\}_{k\in \N}$ in $(L^4(Q_T))^{3\times 3}$ and the convention \eqref{conven}) and the fact that $\bv^1_k$ is a weak solution~$\bu_k$ to the equation \eqref{dva12} with $\bH_k:= g_k (|\bF_k|) \bF_k \bF_k^T$, all the assumptions of Theorem \ref{bitingT} are satisfied. Hence, due to the properties of $E_j$ and the symmetry of $\tilde{\bG}_k$ and $\tilde{\bG}$, we obtain  (for a suitable subsequence of $\{\tilde{\bG}_k:\bD^1_k\}_{k\in \N}$, that we do not relabel), that
\begin{equation}
    \lim_{k\to\infty}\int_{E_j} (\tilde{\bG}_k:\bD^1_k)\varphi= 
    \int_{E_j}(\tilde{\bG}:\bD^1)\varphi.
\end{equation}
Therefore, using again the convergence results \eqref{vk1weak}, \eqref{gkFkstrong} and \eqref{conven} with $\bQ_k:= \bF_k \bF_k^T$, we conclude that
\begin{equation}
\label{L1j1est}
L_1^{j,1}= -\lim_{k\to \infty}\int_{E_j} |\bD_k^1-\bD^1|^2\varphi.
\end{equation}
Next, we  introduce an appropriate estimate for $L_3^{j,1}$. We have, using \eqref{vk1weak}, \eqref{Fkweak}, Young's inequality and the localised form of Korn's inequality (see~\cite{BuMaLo22}, Step 2 of the proof of the compactness of $\{\bF_\varepsilon\}$), that
\begin{equation}
    \label{L3j1est}
    \begin{aligned}
&L_3^{j,1}= \lim_{k\to\infty}\int_{E_j}\left((\nabla\bv^1-\nabla \bv^1_k) (\bF_k-\bF): \bF\right)\varphi
\\&\le \lim_{k\to\infty}\int_{E_j} \left(|\bD^1-\bD^1_k|^2 + |\bF_k-\bF|^2|\bF|^2\right)\varphi\\&\le\lim_{k\to \infty} \int_{E_j}|\bD^1-\bD^1_k|^2 \varphi+\int_{Q_T}|\bF|^2(\overline{|\bF|^2}-|\bF|^2)\varphi.
\end{aligned}
\end{equation}
Note that, summing \eqref{L1j1est} and \eqref{L3j1est}, we obtain, for all $j\in \N$, the estimate
\begin{equation}
    \label{firstmain}
    L_1^{j,1}+L_3^{j,1}\leq \int_{Q_T}|\bF|^2(\overline{|\bF|^2}-|\bF|^2)\varphi.
\end{equation}

Next, we show that
\begin{equation}
    \label{Ltj1}\lim_{j\to\infty}L_{1}^{\tilde{j},1}=\lim_{j\to\infty}L_{3}^{\tilde{j},1}=0.
\end{equation}
To prove the result for $L_3^{\tilde{j},1}$ is simple. Thanks to the uniform boundedness of $\{\nabla \bv_k^1\}_{k\in \N}$ in $(L^2(Q_T))^{3\times3}$ and $\{\bF_k\}_{k\in \N}$ in $(L^4(Q_T))^{3\times 3}$ and the fact that  $\bF\in (L^4(Q_T))^{3\times3}$, we observe  that (use also \eqref{conven})
\begin{equation}
\label{L3t11}
    L_3^{\tilde{j},1} = \int_{Q_T\setminus E_j} (\nabla \bv^1\bF:\bF - \overline{\nabla \bv^1 \bF}:\bF)\varphi\xrightarrow{j\to \infty}0,
\end{equation}
where the last convergence follows from the facts that the integrand on the right-hand side of \eqref{L3t11} belongs to $L^1(Q_T)$ and that $|Q_T\setminus E_j|\to 0$ as $j\to \infty$.

Hence, we can focus on the term $L_1^{\tilde{j},1}$. Similarly as above, 
we observe that
\begin{equation}
    \label{Ltj1pom0}\lim_{k\to \infty}\int_{Q_T\setminus E_j}\left(\bD^1: g_k(|\bF_k|) \bF_k \bF_k^T\right)\varphi=\int_{Q_T\setminus E_j}(\bD^1:\overline{\bF \bF^T})\varphi\xrightarrow{j\to \infty}0,
    \end{equation}
and thus, to prove that $L_1^{\tilde{j},1}\to 0$ as $j\to\infty$, it suffices to show  the identity
\begin{equation}
\label{Ltj1pom}
\lim_{j\to \infty}    \lim_{m\to \infty} \lim_{k\to \infty}\int_{Q_T\setminus E_j} \left(\bD^1_k: g_m(|\bF_k|^2)g_k(|\bF_k|)\bF_k \bF_k^T\right)\varphi=0.
    \end{equation}
Due to the presence of the term $g_m(|\bF_k|^2)$, the uniform boundedness of $\{\bD^1_k\}_{k\in \N}$ in $(L^2(Q_T))^{3\times 3}$ and due to \eqref{gkFkstrong} (and \eqref{conven}), we have the following convergence results as $k\to\infty$:
\begin{align}
\label{DmFFweak}
\bD^1_k: g_m(|\bF_k|^2)g_k(|\bF_k|)\bF_k \bF_k^T&\rightharpoonup \overline{\bD^1: g_m(|\bF|^2)\bF \bF^T} &&\textrm{weakly in }L^2(Q_T),\\
\label{dmnFweak}
|\bD^1_k|g_m(|\bF_k|^2)g_k(|\bF_k|)|\bF_k|^2&\rightharpoonup \overline{|\bD^1| g_m(|\bF|^2)|\bF|^2} &&\textrm{weakly in }L^2(Q_T).
\end{align}
Obviously, by the weak-lower semicontinuity of $L^1$-norm, one has that
\begin{equation}
    \label{FFnFest} \left|\overline{\bD^1:g_m(|\bF|^2)\bF \bF^T}\right|\leq \overline{|\bD^1|g_m(|\bF|^2)|\bF|^2}
    \quad \mbox{ a.e. in } Q_T
\end{equation}and, since $g_m\le g_n$ whenever $m\le n$, we have also that 
\begin{equation}
\overline{|\bD^1| g_m(|\bF|^2)|\bF|^2}\le \overline{|\bD^1| g_n(|\bF|^2)|\bF|^2} \quad \mbox{a.e. in } Q_T.\label{TL2}
\end{equation} 
As the sequence $\{\overline{|\bD^1|g_m(|\bF|^2)|\bF|^2}\}_{m\in \N}$ is uniformly bounded in $L^1(Q_T)$, we can employ \eqref{TL2} and the monotone convergence theorem  to conclude that
there exists $\overline{|\bD^1|g_\infty(|\bF|^2)|\bF|^2}$ such that the following result holds as $m\to \infty$:
\begin{equation}
    \label{TL3}
    \overline{|\bD^1| g_m (|\bF|^2)|\bF|^2}\to \overline{|\bD^1| g_\infty(|\bF|^2)|\bF|^2}\quad \mbox{strongly in } L^1(Q_T).
\end{equation}
This means that the sequence $\{\overline{|\bD^1| g_m (|\bF|^2)|\bF|^2}\}_{m\in \N}$ is uniformly equi-integrable on $Q_T$. Hence, by \eqref{FFnFest}, the sequence $\{\overline{\bD^1:g_m(|\bF|^2)\bF\bF^T}\}_{m\in \N}$ is also uniformly equi-integrable on $Q_T$. This fact together with the property $|Q_T\setminus E_j|\to 0$ as $j\to \infty$ and with \eqref{DmFFweak} leads to the desired identity \eqref{Ltj1pom}, which, together with \eqref{L3t11} and \eqref{Ltj1pom0}, implies the validity of \eqref{Ltj1}. 

Last, we prove that
\begin{equation}
\label{Lv2}
    L^2_1=L^2_3=0.
\end{equation} The result $L_3^2=0$ follows immediately from \eqref{vk2strong}, \eqref{vk2weak} and \eqref{Fkweak}. Thus, it suffices to prove that $L_1^2=0$. We have that
\begin{equation*}
\begin{aligned}
&\left|\lim_{m\to \infty} \lim_{k\to\infty} \int_{Q_T} \left(\bD^2: (g_m(|\bF_k|^2)-1) g_k(|\bF_k|) \bF_k \bF_k^T\right)\varphi\right|\\ &\qquad\leq \|\varphi\|_{L^\infty(Q_T)}\lim_{m\to \infty} \lim_{k\to\infty} \int_{\{|\bF_k|^2>m\}} |\bF_k|^2|\bD^2|=0,
\end{aligned}
\end{equation*}
where the equality holds thanks to \eqref{measFm}, the fact that $\bD^2\in (L^2(Q_T))^{3\times 3}$, the uniform boundedness of $\{\bF_k\}_{k\in \N}$ in $(L^4(Q_T))^{3\times 3}$ and Hölder's inequality. Hence
\begin{equation}
\label{L12equiv}
    L_1^2= \lim_{m\to \infty} \lim_{k\to\infty} \int_{Q_T} \left((g_m(|\bF_k|^2) g_k(|\bF_k|) \bF_k \bF_k^T):(\bD^2_k-\bD^2)\right)\varphi.
\end{equation}
Thanks to the presence of the cut-off function $g_m$ and the strong convergence \eqref{vk2strong}, we conclude from \eqref{L12equiv} that
\begin{equation}
    \label{L12concl}
    |L_1^2|\leq \lim_{m\to\infty}\left(2m \lim_{k\to \infty} \int_{Q_T} |\bD^2_k-\bD^2| \varphi\right) = 0. 
\end{equation} Hence (recall that $L_3^2=0$) the relation \eqref{Lv2} holds. 

\newpage
To conclude, using the decomposition \eqref{L13dec} and summing the identities \eqref{firstmain}, \eqref{Ltj1} and \eqref{Lv2}, we get the estimate
$$
L_1+L_3\le \int_{Q_T}|\bF|^2 (\overline{|\bF|^2}-|\bF|^2) \varphi,
$$
which, combined with \eqref{step2Fsigmaprep}, yields that \eqref{step2compFsigma} holds with
$$
\tilde{L}:=|\bD|+|\bF|^2 \in L^2(Q_T).
$$
\end{proof}

\begin{proof}[\bf Ad Step 3]
As a result of the previous steps, we have that
\begin{equation}
    -\int_{Q_T} (\overline{|\bF|^2}-|\bF|^2)\partial_t \varphi-\int_{Q_T} (\overline{|\bF|^2}-|\bF|^2)\bv\cdot \nabla \varphi\leq \int_{Q_T} (1+2|\bD|+2|\bF|^2) (\overline{|\bF|^2}-|\bF|^2) \varphi.
\end{equation}
Hence, employing Proposition~\ref{convae} with $f:= \overline{|\bF|^2}-|\bF|^2\ge0 $, $\bu:=\bv$ and $L:=(1+2|\bD|+2|\bF|^2)$, we get~\eqref{compae}. The relation \eqref{compae} is equivalent to the compactness $\{\bF_k\}_{k\in \N}$ in $(L^2(Q_T))^{3\times3}$.
Indeed,
$$
\lim_{k\to\infty} \|\bF_k-\bF\|_{2,Q_T}^2 = \lim_{k\to\infty} \int_{Q_T} \left(|\bF_k|^2 - 2 \bF_k: \bF + |\bF|^2\right)= \int_{Q_T} \left(\overline{|\bF|^2} - |\bF|^2\right)=0,
$$
where the last equality follows from \eqref{compae}. The proof of the compactness of
$\{\bF_k\}_{k\in \N}$ in $(L^2(Q_T))^{3\times3}$ is complete.
\end{proof}

\subsection{Positivity of \texorpdfstring{$\det \bF$}{detF}} \label{limitk2}
As we have proved the compactness of $\{\bF_k\}_{k\in \N}$ in $(L^2(Q_T))^{3\times3}$, we can find a (not relabeled) subsequence of $\{\bF_k\}_{k\in \N}$ such that $\bF_k\to \bF$ a.e. in $Q_T$. Then from Fatou's lemma and the estimate~\eqref{3sta} we conclude that
\begin{equation}\label{3stae}
\begin{split}\sup_{t\in (0,T)}\|\ln \det \bF(t)\|_1\le\sup_{t\in (0,T)}\liminf_{k\to\infty}\|\ln \det \bF_k(t)\|_1 \le C(T, \Omega,\bF_0,\bv_0),
\end{split}
\end{equation}
which gives that $\det \bF>0$ a.e. in~$Q_T$. The proof of Theorem \ref{main} is complete.

%

\section{Extensions}\label{Extensions}
In this section, we generalize Theorem~\ref{Burgmain} (the basic existence result for the classical Giesekus model in three dimensions) in several directions. First, we cover all spatial dimensions $d\ge 2$. Second, we consider models with two relaxation mechanisms developing thus a mathematical theory for generalizations of the second order viscoelastic models such as the classical model due to Burgers. Third, we also consider a whole hierarchy of rate-type fluid models that include the Giesekus model as a limiting case and that, as will be shown below, can be analyzed in a similar fashion. All these three types of generalization are considered  together in what follows.

Let $\Omega$ and $T$ satisfy \eqref{p000}. Given the parameters $G_1$, $G_2$, $\tau_1$, $\tau_2$, $\rho$, $\nu\in (0,\infty)$, $\lambda_1$, $\lambda_2\in\mathbb{R}$ and the external body forces $\vec{f}:Q_T\to \R^d$, we look for functions $(\bv, p, \bF_1, \bF_2, \bB_1, \bB_2): Q_T\to \R^d\times \R \times \R^{d\times d} \times \R^{d\times d} \times \R^{d\times d} \times \R^{d\times d}$ satisfying the following governing equations and the natural constraints in $Q_T$:
\small\begin{align}
\label{1Burgi}
\diver \bv&=0,\\
\rho\partial_t \bv + \rho\diver (\bv\otimes \bv) + \nabla p - \diver 2\nu\bD+\sum_{i=1}^2 G_i \diver \bB_i &=\rho\vec{f},\\
\label{vB3lambdai}
\partial_t \bB_i +\mbox{Div} (\bB_i\otimes \bv) - \nabla \bv \bB_i -\bB_i (\nabla \bv)^T+ \frac{1}{\tau_i}\left(\bB_i^{2-\lambda_i} - \bB_i^{1-\lambda_i}\right)&=\tens{O},&& i=1,2,\\
\bF_i \bF_i^T&=\bB_i,&& i=1,2, \\
\det \bF_i&>0,&& i=1,2,  \label{detFformi}
\end{align}
\normalsize together with the homogeneous Dirichlet boundary condition for $\bv$:
\begin{equation}\label{Burgbci}
\bv=\vec{0} \quad \mbox{on } (0,T)\times \partial \Omega,
\end{equation}
and the following initial conditions for $\bv$, $\bF_i$ and $\bB_i$ ($i=1,2$):
\begin{equation}\label{ici}
\bv(0,\vec{\cdot})=\bv_0, \quad \bF_i(0,\vec{\cdot})=\bF_{i_0}, \quad \bB_i(0,\vec{\cdot})=\bB_{i_0}:= \bF_{i_0} \bF_{i_0}^T \quad \mbox{in } \Omega.
\end{equation}

If $\lambda_1=\lambda_2=1$, then the system \eqref{1Burgi}--\eqref{detFformi} becomes the classical Burgers viscoelastic rate-type fluid model of second order written as a ``mixture" (a combination) of two classical Oldroyd-B viscoelastic rate-type fluid models of the first order, see e.g. \cite{MRT} for a proof of this property.  

Let us recall that rate-type viscoelastic fluid models of second order (such as the~Burgers model and its generalizations) are the simplest models that are capable of describing responses of complex geo-, bio- and polymeric materials, see e.g. \cite{Reh,Nar,derivold,Sk,Tuma18}.

It is also interesting to notice under what assumption we obtain the classical Giesekus model studied in previous sections. It requires setting $G_2 =0$ (implying that the evolutionary equation for $\bB_2$ is redundant), setting $\lambda=\lambda _1 =0$ and relabelling $G_1$ and $\tau_1$ to $G$ and $\tau$. 
As one can see bellow, we will be able to extend the developed approach just to all $\lambda\le 0$ (the method is not applicable to the Oldroyd-B model, i.e. $\lambda = 1$). This is the reason why, in what follows, we set 
$\mu_i:=-\lambda_i$, $i=1,2$. 
Let $\vec{f}=\vec{0}$ for simplicity. The system \eqref{1Burgi}--\eqref{ici} satisfies, for all $t>0$, the following a~priori energy identity (supposed that all the integrated quantities are smooth enough, see \cite[Relation (2.20)]{BuMaLo22} for details):
\begin{equation} 
\label{Pepa41}
\begin{split}
&\frac12\int_{\Omega} \left(\rho|\bv(t)|^2 + \sum_{i=1}^2 G_i( \tr \bB_{i}(t) -d-\ln \det \bB_{i}(t))\right) \\
&\qquad +  \int_0^t \int_{\Omega} \left(2 \nu |\bD|^2 + \sum_{i=1}^2  \frac{G_i}{2\tau_i} |\bB_{i}^{\frac{2+\mu_i}{2}}(\bI-\bB^{-1}_{i})|^2\right)\\
&=\frac12\int_{\Omega} \left( \rho|\bv_0|^2 + \sum_{i=1}^2 G_i (\tr \bB_{i_0} - d- \ln \det \bB_{i_0})\right).
\end{split}
\end{equation}
It is natural to require that the initial data are such that the right-hand side of~\eqref{Pepa41} is finite. To guarantee that all of the terms appearing in the weak formulation of \eqref{1Burgi}--\eqref{vB3lambdai} are integrable, one observes that $\nabla \bv \bB_i$, resp. $\bB_i (\nabla \bv)^T$, $i=1,2$, are the most critical terms. By \eqref{Pepa41} and Korn's inequality we may deduce that $\nabla \bv\in (L^2(Q_T))^{d\times d}$ and $\bB_i\in (L^{2+\mu_i}(Q_T))^{d\times d}$. Hence, to conclude  that $\nabla \bv \bB_i$ and $\bB_i (\nabla \bv)^T$ belong to $(L^1(Q_T))^{d\times d}$, we need that $\bB_i\in (L^2(Q_T))^{d\times d}$, which is true if $\mu_i\geq 0$.
We aim to prove the existence of a weak solution to the problem \eqref{1Burgi}--\eqref{ici} for any $\mu_1, \mu_2\in [0, \infty)$. The result is stated in the following theorem.
\begin{thm}
\label{Burgmaini}
Let $\Omega$ and $T$ satisfy \eqref{p000} and $G_1, G_2, \tau_1, \tau_2, \rho, \nu\in(0,\infty)$ be given constants. Let $\mu_1, \mu_2 \in [0, \infty)$. Let $\vec{f}\in L^2(0,T;(W_{\vec{0},\diver}^{1,2})^*)$, $\bv_0\in L^2_{\bn, \diver}$, $\bF_{i_0}\in (L^2(\Omega))^{d\times d}$, $\det \bF_{i_0}>0$ a.e. in $\Omega$ and $\ln \det \bF_{i_0}\in L^1(\Omega)$ ($i=1,2$).  Then, there exists a weak solution to \eqref{1Burgi}--\eqref{ici}, i.e. there exists a quintuple $(\bv, \bF_1, \bF_2,  \bB_1, \bB_2)$ satisfying (for $i=1,2$)
\begin{align*}
\bv & \in \mathcal{C}_{weak}([0,T]; L^2_{\bn, \diver}) \cap L^2(0,T; W_{\vec{0},\diver}^{1,2}), \\
\partial_t \bv &\in L^{\frac{d+2}{d}}\left(0, T; (W_{\vec{0},\diver}^{1,\frac{d+2}{2}})^*\right),\\
\bB_i & \in \mathcal{C}([0,T]; (L^1(\Omega))^{d\times d}) \cap (L^{2+\mu_i}(Q_T))^{d\times d},
\\
\partial_t \bB_i & \in L^1 (0,T; ((W^{1,d+1}(\Omega))^{d\times d})^*),
\\
\bF_i \bF_i^T&= \bB_i \mbox{ a.e. in } Q_T,
\\
\det \bF_i &> 0 \mbox{ a.e. in } Q_T,
\end{align*}
such that, for all $\bw \in W_{\vec{0},\diver}^{1,\frac{d+2}{2}}$, all $\bA \in (W^{1,d+1}(\Omega))^{d\times d}$ and a.a. $t\in (0,T)$, the following identities hold:
\small\begin{align}
\label{vBweaklambda1i}
&\rho\langle \partial_t \bv, \bw \rangle- \rho\int_{\Omega} (\bv \otimes \bv): \nabla \bw + \int_{\Omega} \left(2\nu\bD+\sum_{i=1}^2G_i\bB_i\right): \nabla \bw-\langle \vec{f}, \bw\rangle= 0,\\
\label{vBweaklambda2i}
&\langle \partial_t \bB_i, \bA \rangle -\int_{\Omega} (\bB_i \otimes \bv): \nabla \bA +  \left(\nabla \bv \bB_i+\bB_i (\nabla \bv)^T - \frac{1}{\tau_i} (\bB_i^{2+\mu_i} -  \bB_i^{1+\mu_i})\right): \bA = 0,
\end{align}
\normalsize and the initial data $\bv_0$, $\bF_{i_0}$, $\bB_{i_0}:= \bF_{i_0} \bF_{i_0}^T$ are attained in the following sense:
\begin{equation}
\label{invFBweaki}
\lim_{t\to 0+} \left(\|\bv(t)-\bv_0\|_2 +\sum_{i=1}^2\|\bF_i(t)-\bF_{i_0}\|_2 + \|\bB_i(t)-\bB_{i_0}\|_1\right)= 0.
\end{equation}
\end{thm}

Let us note that, having $\bv\in L^\infty(0,T; L^2_{\bn,\diver})\cap L^2(0,T; W_{\vec{0},\diver}^{1,2})$ and using the standard interpolation and the Sobolev embedding, we obtain that \begin{equation}\label{interpolv}\bv\in (L^\frac{2(d+2)}{d}(Q_T))^d,\end{equation} hence, due to the Hölder inequality, the function $(\bv\otimes \bv):\nabla \bw$ is integrable over $Q_T$ for arbitrary $\bw\in L^{(d+2)/2}(0,T; W_{\vec{0}, \diver}^{1, (d+2)/2})$. Since this space is embedded to $ L^2(0,T; W_{\vec{0},\diver}^{1,2})$, all the terms acting in \eqref{vBweaklambda1i} are well defined and integrable over $(0,T)$.\footnote{As one may check, if $d>4$, then it is not possible to ensure that $\partial_t \bv\in L^1(0,T;(W_{\vec{0},\diver}^{1,2})^*)$ and test \eqref{vBweaklambda1i} by $\bw\in W_{\vec{0},\diver}^{1,2}$ for a fixed $t\in (0,T)$. The function space for $\partial_t \bv$ is chosen so that $\nabla \bw$ in \eqref{vBweaklambda1i} integrated over $(0,T)$ has the same integrability properties in time and space.}  Next, having $\bB_i\in (L^{2+\mu_i}(Q_T))^{d\times d}$, we need the test function~$\bA$ in the equation \eqref{vBweaklambda2i} integrated over $(0,T)$ to be in $(L^\infty(Q_T))^{d\times d}$. Due to the properties $\bv\in L^2(0,T; (L^\frac{2d}{d-2}(\Omega))^d)$ (valid as $W^{1,2}(\Omega)\hookrightarrow L^\frac{2d}{d-2}(\Omega)$), $\bB_i\in (L^2(Q_T))^{d\times d}$ and the Hölder inequality the term $(\bB_i\otimes \bv):\nabla \bA$ is integrable for any $\bA\in L^\infty(0,T;(W^{1,d}(\Omega))^{d\times d})$. Since $W^{1, d+1}(\Omega)$ is embedded to $ L^\infty (\Omega)$ we may set $\bA \in L^\infty(0,T;(W^{1,d+1}(\Omega))^{d\times d})$ in the equation \eqref{vBweaklambda2i} integrated over $(0,T)$. This is why $\partial_t \bB_i$ shall belong to $L^1(0,T; ((W^{1,d+1}(\Omega))^{d\times d})^*)$ and all the terms acting in \eqref{vBweaklambda2i} are well defined and integrable over $(0,T)$. 

In the rest of the paper we focus on the proof of Theorem \ref{Burgmaini}. As in the proof of Theorem \ref{Burgmain}, we set for simplicity and without loss of any generality the material constants $\rho$, $2\nu$, $\tau_1$, $\tau_2$,  $G_1$, $G_2$ to be equal to one and $\vec{f}=\vec{0}$. Also, since the  argument how the existence of a weak solution to the mixture of two Giesekus models ($\mu_1, \mu_2=0$ in \eqref{vB3lambdai}) follows from the existence result for the single Giesekus model 
in two spatial dimensions was performed in \cite{BuMaLo22} in detail and the procedure is independent of dimension and also of the value of the parameters $\mu_1$, $\mu_2$, we skip this reasoning here and we focus just on the proof of the existence of a weak solution for a single fluid model with $\mu:=\mu_1\in [0, \infty)$ arbitrary. To put it differently, our objective is to prove the existence of $(\bv, \bF:= \bF_1, \bB:= \bB_1)$ satisfying \eqref{vBweaklambda1i}--\eqref{invFBweaki} with $\rho=2\nu=\tau_1=G_1=1$, $G_2=0$, $\vec{f}=\vec{0}$ and $\mu:=\mu_1\in [0,\infty)$.

Similarly as in the proof of Theorem \ref{Burgmain},  instead of the equation for $\bB$ we consider the evolution equation for $\bF$ that takes, for 
$\mu\in[0,\infty)$, the form
\begin{equation}
    \label{Flambda}
    \partial_t \bF + \mbox{Div } (\bF\otimes \bv) - \nabla \bv \bF + \frac{1}{2} \left((\bF \bF^T)^{1+\mu}-(\bF \bF^T)^{\mu}\right)\bF=\tens{O}.
\end{equation}
Multiplying equation~\eqref{Flambda} by $\bF^T$ from the~right and the transpose of the equation~\eqref{Flambda} by $\bF$ from the~left and summing the results, we arrive at \eqref{vB3lambdai} with $i=1$, $\bB_1=\bB=\bF \bF^T$, $\tau=\tau_1=1$ and $\mu=\mu_1$. This apparently informal procedure can be done rigorously following step by step the corresponding procedure for $\mu:=\mu_1=0$ presented in \cite{BuMaLo22}. Hence the proof of Theorem \ref{Burgmaini} is reduced to the proof of the existence of a weak solution $(\bv, \bF)$ to the problem \eqref{vF1}--\eqref{vF2}, \eqref{Flambda}, \eqref{detF}--\eqref{icvF}.
\begin{thm}
\label{maini}
Let $\Omega$ and $T$ satisfy \eqref{p000} and $\mu\in[0, \infty)$.  Let $\bv_0\in L^2_{\bn, \diver}$ and let $\bF_0\in (L^2(\Omega))^{d\times d}$ be such that $\det \bF_0>0$ a.e. in $\Omega$ and $\ln \det \bF_0\in L^1(\Omega)$. Then, there exists $(\bv, \bF)$, a weak solution to \eqref{vF1}--\eqref{vF2}, \eqref{Flambda}, \eqref{detF}--\eqref{icvF}, satisfying
\begin{align*}
\bv & \in \mathcal{C}_{weak}([0,T]; L^2_{\bn, \diver}) \cap L^2(0,T; W_{\vec{0},\diver}^{1,2}), \\
\partial_t \bv &\in L^{\frac{d+2}{d}}\left(0, T; (W_{\vec{0},\diver}^{1,\frac{d+2}{2}})^*\right),\\
\bF & \in \mathcal{C}([0,T]; (L^2(\Omega))^{d\times d}) \cap (L^{4+2\mu}(Q_T))^{d\times d}, \\
\partial_t \bF & \in L^{\frac{4+2\mu}{3+2\mu}}\left(0,T; ((W^{1,d}(\Omega))^{d\times d})^*\right),\\
\det \bF &> 0 \mbox{ a.e. in } Q_T,\\
\ln (\det \bF) & \in L^{\infty}(0,T; L^1(\Omega)),
\end{align*}
such that, for all $\bw \in W_{\vec{0},\diver}^{1,\frac{d+2}{2}}$, $\bA \in (W^{1,d}(\Omega))^{d\times d}$ and a.a. $t\in (0,T)$, the following identities hold:
\begin{align}
\label{vweaki}
&\langle \partial_t \bv, \bw \rangle- \int_{\Omega} (\bv \otimes \bv): \nabla \bw + \int_{\Omega} (\bD+\bF \bF^T): \nabla \bw= 0,\\
&\begin{aligned}\label{Fweaki}
&\langle \partial_t \bF, \bA \rangle -\int_{\Omega} (\bF \otimes \bv): \nabla \bA - \int_{\Omega} \nabla \bv \bF: \bA \\&\qquad \quad + \frac{1}{2}\int_{\Omega} \left((\bF \bF^T)^{1+\mu}  -  (\bF\bF^T)^{\mu}\right)\bF: \bA = 0,
\end{aligned}
\end{align}
and the initial conditions $\bv_0$ and $\bF_{i_0}$ are attained in the following sense
\begin{equation}
\label{invFweaki}
\lim_{t\to 0+} \left(\|\bv(t)-\bv_0\|_2 + \|\bF(t)-\bF_0\|_2\right) = 0.
\end{equation}
\end{thm}

The rest of this paper is devoted to the proof of Theorem \ref{maini}. A crucial task is to prove the compactness of the approximations to $\bF$ in $(L^2(Q_T))^{d\times d}$. 

\subsection{Mathematical tools}
To be able to treat the terms $(\bF \bF^T)^{1+\mu}\bF$ and $(\bF \bF^T)^{\mu}\bF$ in \eqref{Fweaki} (the terms that make the evolutionary equation for $\bF$ different from the Giesekus model), we establish  the following proposition on their monotonicity and interpolation.

\begin{prop}[Properties of $(\bA \bA^T)^s\bA$]\label{mon}
Let $s\in (-\frac12,\infty)$ be fixed, but arbitrary. Then the matrix function
$S_s(\bA):= (\bA \bA^T)^s \bA$ is monotone:  
\begin{equation}
    \label{monotfin1}
\textrm{for all } \bA,\bC\in \R^{d\times d}: \quad (S_s(\bA)-S_s(\bC)):(\bA-\bC)\geq 0.
\end{equation}
Moreover, if $s\ge 0$, then, for any $\varepsilon>0$, there is a constant $c=c(s,\varepsilon)\in (0,\infty)$ such that
\begin{equation}
\label{monotfin2}(S_s(\bA)-S_s(\bC)):(\bA-\bC)\leq \varepsilon(S_{s+1}(\bA)-S_{s+1}(\bC)):(\bA-\bC)
+c |\bA-\bC|^2.
\end{equation}
\end{prop}
\begin{proof}
Let $\bA, \bC\in \R^{d\times d}$ be arbitrary. Let 
\begin{equation}
\label{monotdefAC}
    \bA=\tens{U} \tens{H} \tens{V}^T,\quad 
  \bC=\tilde{\tens{U}} \tilde{\tens{H}}\tilde{\tens{V}}^T
\end{equation} 
be  their singular value decompositions with $\tens{U}$, $\tens{V}$, $\tilde{\tens{U}}$, $\tilde{\tens{V}}$ being orthogonal matrices (matrices $\tens{X}\in \R^{d \times d}$ satisfying $\tens{X} \tens{X}^T=\bI$) and $\tens{H}$, $\tilde{\tens{H}}$ being of the form
\begin{equation}
\label{monotdefH}
(\bH)_{i,j}:= \delta_{ij} \sigma_j, \quad (\tilde{\bH})_{i,j}:= \delta_{ij} \tilde{\sigma}_j,\quad i,j=1,...,d,
\end{equation}
where $\delta_{ij}$ stands for the Kronecker symbol and $\sigma_j:=\sqrt{\kappa_j}$ and $\tilde{\sigma}_j:=\sqrt{\tilde{\kappa}_j}$ ($j=1,...,d$) are the singular values of $\bA$ and $\bC$, respectively (which means that $\kappa_j$ and $\tilde{\kappa}_j$, $j=1,...,d$, are the eigenvalues of $\bA \bA^T$ and $\bC \bC^T$, respectively). 
From \eqref{monotdefAC}, \eqref{monotdefH} and the orthogonality of $\bU$, $\tilde{\bU}$, $\bV$ and $\tilde{\bV}$ it follows that 
\begin{equation}\label{AAT} (\bA \bA^T)^s= \tens{U} \tens{H}^{2s} \tens{U}^T, \quad (\bC\bC^T)^s= \tilde{\tens{U}} \tilde{\tens{H}}^{2s} \tilde{\tens{U}}^T,
\end{equation}
\begin{equation}\label{AAT2} S_s(\bA)= \tens{U} \tens{H}^{2s+1} \tens{V}^T, \quad S_s(\bC)= \tilde{\tens{U}} \tilde{\tens{H}}^{2s+1} \tilde{\tens{V}}^T.
\end{equation}
Note that \eqref{AAT} is well defined only if all singular values $\sigma_j$, $\tilde{\sigma}_j$ are strictly positive or if $s\ge 0$, while \eqref{AAT2} is well defined whenever $s>-\frac12$.
Moreover, the orthogonality of $\bU$, $\bV$, $\tilde{\bU}$ and $\tilde{\bV}$ implies the orthogonality of the matrices $\tilde{\bU}^T \bU$ and $\tilde{\bV}^T \bV$, and thus, for all $i,j=1,...,d$, it holds
\small\begin{equation}\label{sigmaij} 1= \frac{1}{2}\sum_{l=1}^d\left((\tilde{\bU}^T \bU)_{l,j}^2+ (\tilde{\bV}^T \bV)_{l,j}^2\right)= \frac{1}{2}\sum_{l=1}^d\left((\tilde{\bU}^T \bU)_{i,l}^2+ (\tilde{\bV}^T \bV)_{i,l}^2\right). 
\end{equation}
\normalsize  Using \eqref{monotdefAC}--\eqref{AAT2} and the orthogonality of $\bU$, $\bV$, $\tilde{\bU}$ and $\tilde{\bV}$, we have
\small\begin{align*}
    &(S_s(\bA)-S_s(\bC)):(\bA-\bC)=S_s(\bA):\bA+S_s(\bC):\bC-S_s(\bA):\bC-S_s(\bC):\bA
    \\&= \bH^{2s+1}:\bH+\tilde{\bH}^{2s+1}:\tilde{\bH}-(\tilde{\bU}^T \bU) \bH^{2s+1}:\tilde{\bH}(\tilde{\bV}^T \bV)-(\tilde{\bU}^T \bU)\bH:\tilde{\bH}^{2s+1}(\tilde{\bV}^T \bV)
    \\ &=\sum_{j=1}^d(\sigma_j^{2s+2}+\tilde{\sigma}_j^{2s+2})-\sum_{i,j=1}^d (\sigma_j^{2s+1}\tilde{\sigma}_i+\sigma_j\tilde{\sigma}_i^{2s+1})(\tilde{\bU}^T \bU)_{i,j}(\tilde{\bV}^T \bV)_{i,j}.
\end{align*}
\normalsize Employing \eqref{sigmaij}, we see that the last expression is equal to
\small
\begin{align*}
&\frac{1}{2}\sum_{i,j=1}^d (\sigma_j^{2s+2}+\tilde{\sigma}_i^{2s+2}-\sigma_j^{2s+1}\tilde{\sigma}_i-\sigma_j\tilde{\sigma}_i^{2s+1})\left((\tilde{\bU}^T
\bU)^2_{i,j}+(\tilde{\bV}^T \bV)^2_{i,j}\right)
\\ & \quad +\frac{1}{2}\sum_{i,j=1}^{d}(\sigma_j^{2s+1}\tilde{\sigma}_i+\sigma_j \tilde{\sigma}_i^{2s+1})\left((\tilde{\bU}^T\bU)_{i,j}-(\tilde{\bV}^T \bV)_{i,j}\right)^2.
\end{align*}
\normalsize
To summarize, the above computation leads to 
\small\begin{equation}\label{monotkey}
\begin{split}
(S_s(\bA) &-S_s(\bC)):(\bA-\bC) \\ &=\frac{1}{2}\sum_{i,j=1}^d(\sigma_j^{2s+1}-\tilde{\sigma}_i^{2s+1})(\sigma_j-\tilde{\sigma}_i)\left((\tilde{\bU}^T \bU)_{i,j}^2+(\tilde{\bV}^T \bV)_{i,j}^2\right)\\ & \qquad +(\sigma_j^{2s+1}\tilde{\sigma}_i+\sigma_j \tilde{\sigma}_i^{2s+1})\left((\tilde{\bU}^T\bU)_{i,j}-(\tilde{\bV}^T \bV)_{i,j}\right)^2.\end{split}
\end{equation}
\normalsize
Since $g(y):= y^{2s+1}$ is increasing on $[0,\infty)$ and the singular values $\sigma_j$ and $\tilde{\sigma}_j$, $j=1,...,d$, are nonnegative, each summand on the right-hand side of \eqref{monotkey} is nonnegative. This completes the proof of  \eqref{monotfin1}. 

To prove \eqref{monotfin2}, it suffices to show, taking into account \eqref{monotkey} and the nonnegativity of each summand on its right-hand side, that ($i,j=1,...,d$ and $c=c(s,\varepsilon)$)
\begin{align}
\label{monotfin2a}
(\sigma_j^{2s+1}-\tilde{\sigma}_i^{2s+1})(\sigma_j-\tilde{\sigma}_i)&\leq \varepsilon(\sigma_j^{2s+3}-\tilde{\sigma}_i^{2s+3})(\sigma_j-\tilde{\sigma}_i)+c(\sigma_j-\tilde{\sigma}_i)^2,
\\
\label{monotfin2b}
(\sigma_j^{2s+1}\tilde{\sigma}_i+\sigma_j \tilde{\sigma}_i^{2s+1})&\leq \varepsilon(\sigma_j^{2s+3}\tilde{\sigma}_i+\sigma_j \tilde{\sigma}_i^{2s+3})+c(\sigma_j\tilde{\sigma}_i+\sigma_j \tilde{\sigma}_i).
\end{align}
The inequality \eqref{monotfin2b} with $c(s,\varepsilon):= \varepsilon^{-s}$ follows from the nonnegativity of each $\sigma_j$, $\tilde{\sigma}_j$ ($j=1,...,d$) and from the fact that 
$$
y^{2s}\leq \varepsilon y^{2s+2} +\varepsilon^{-s}
$$ 
valid for all $y\ge 0$.

It remains to prove \eqref{monotfin2a}. We can assume, without loss of any generality, that $\tilde{\sigma}_i\leq \sigma_j$. 
  Thus, to prove \eqref{monotfin2a}, we need to show that $f(t)$ defined as 
$$
f(t):= \varepsilon (t^{2s+3}-t_0^{2s+3})-(t^{2s+1}-t_0^{2s+1}) + c(s,\varepsilon)(t-t_0)
$$
is nonnegative for $t\in [t_0,\infty)$. Since $f(t_0)=0$, it is enough to show that $f$ is nondecreasing. Thus, applying the derivative, we have that
$$
f'(t)=\varepsilon(2s+3) t^{2s+2}-(2s+1)t^{2s}+ c(s,\varepsilon).
$$
We see that if 
\begin{equation}
\label{ces}    
c(\varepsilon,s):=(2s+1)\left(\frac{2s+1}{\varepsilon(2s+3)}\right)^s,
\end{equation}then $f'$ is nonnegative for any $t\ge t_0$. Hence, the proof of \eqref{monotfin2a} is complete.


To conclude, \eqref{monotfin2a} holds 
with  $c=c(s,\varepsilon)$ defined by \eqref{ces}, which, together with
\eqref{monotfin2b}, where $c=\varepsilon^{-s}$, and \eqref{monotkey}, gives \eqref{monotfin2} with $c$ equal to the maximum of $\varepsilon^{-s}$ and the value defined in \eqref{ces}.
\end{proof}

\subsection{Proof of Theorem \ref{maini} (main ideas)}

\subsubsection{\texorpdfstring{$(\varepsilon,k)$}{ejk}-approximations and the limits \texorpdfstring{$\varepsilon\to 0+$}{epsk} and \texorpdfstring{$k\to \infty$}{ep}}

\noindent 
For any $\varepsilon \in(0,1)$ and $k\in \mathbb{N}$, we look for $(\bv, p, \bF): Q_T\to \mathbb{R}^d\times \mathbb{R}\times\mathbb{R}^{d\times d}$ satisfying 
\small
\begin{align}
\label{approxepssigma1i}
\diver \bv & = 0,\\
\label{approxepssigma2i}
\partial_t \bv + \diver(g_k(|\bv|) \bv\otimes \bv)+ \nabla p - \diver(\bD + g_k(|\bF|)\bF\bF^T)& = \vec{0},\\
\label{approxepssigma3i}
\partial_t \bF + \mbox{Div } (\bF\otimes \bv) - g_k (|\bF|) (\nabla \bv) \bF + \frac{1}{2} \left((\bF \bF^T)^{1+\mu}  - (\bF \bF^T)^{\mu}\right)\bF & = \varepsilon \Delta \bF,
\end{align}
together with the boundary and initial conditions
\begin{align}
    \label{bcapproxepssigmai}
    &\bv= \vec{0} \quad\mbox{and}  \quad \nabla \bF\cdot \bn:= (\nabla F_{i,j}\cdot \bn)_{i,j=1}^d = \bO \quad\mbox{on } (0,T)\times \partial\Omega, \\
    \label{icapproxepssigmai} &\bv(0, \vec{\cdot}) = \bv_0 \quad \mbox{and} \quad \bF(0,\vec{\cdot}) =\bF_0^k:= \bF_0\chi_{\{|\bF_0|\le k\}}+ \bI\chi_{\{|\bF_0|>k\}} \quad \mbox{in } \Omega.
\end{align}
\normalsize 

Due to the presence of the stress diffusion term on the right-hand side of \eqref{approxepssigma3i} the existence of a weak solution to \eqref{approxepssigma1i}--\eqref{icapproxepssigmai} can be proved by the methods developed e.g. in \cite{BaBuMa21,BuMaPrSu21,BuFeMa19,KrPoSa15}, but also in \cite[Appendix B]{BuMaLo22}. 
 
We shall take the limit as $\varepsilon\to 0+$, then as $k\to \infty$ in \eqref{approxepssigma1i}--\eqref{icapproxepssigmai}. To take the limit as $\varepsilon\to0+$ is simpler (due to the presence of the function $g_k$ for fixed $k\in \N$) and one can follow step by step the limiting procedure as $k\to\infty$ presented below. Note that after letting $\varepsilon\to0+$ the stress diffusion term $\varepsilon \Delta \bF$ on the right-hand side of \eqref{approxepssigma3i} vanishes (and also the boundary condition $\nabla \bF\cdot\bn=\bO$ becomes irrelevant) and the equation for $\bF$ becomes
 \begin{equation}
 \label{approxepssigma3ki}
     \partial_t \bF + \mbox{Div } (\bF\otimes \bv) - g_k (|\bF|)(\nabla \bv) \bF + \frac{1}{2} \left((\bF \bF^T)^{1+\mu}  - (\bF \bF^T)^{\mu}\right)\bF  = \bO.
 \end{equation}
Consequently, we have a weak solution  to \eqref{approxepssigma1i}--\eqref{approxepssigma2i}, \eqref{approxepssigma3ki}, \eqref{Burgbci}, \eqref{icapproxepssigmai}, i.e., there is a couple $(\bv_k, \bF_k)$ such that
\begin{align*}
\bv_k & \in \mathcal{C}([0,T]; L^2_{\bn, \diver}) \cap L^2(0,T; W_{\vec{0},\diver}^{1,2}), \\
\partial_t \bv_k &\in L^{2}(0, T; (W_{\vec{0},\diver}^{1,2})^*),\\
\bF_k & \in \mathcal{C}([0,T]; (L^2(\Omega))^{d\times d}) \cap (L^{4+2\mu}(Q_T))^{d\times d}, \\
\partial_t \bF_k & \in L^{\frac{4+2\mu}{3+2\mu}}(0,T; ((W^{1,d}(\Omega))^{d\times d})^*),
\end{align*}
satisfying, for all $\bw \in W_{\vec{0},\diver}^{1,2}$, $\bA \in (W^{1,d}(\Omega))^{d\times d}$ and a.a. $t\in (0,T)$, the following identities ($\bD_k:=\tfrac12(\nabla \bv_k+(\nabla \bv_k)^T)$):
\begin{align}
\label{Mvweaki}
&\begin{aligned}
\langle \partial_t \bv_k, \bw \rangle &- \int_{\Omega} (g_k(|\bv_k|)\bv_k \otimes \bv_k ) : \nabla \bw \\
&+ \int_{\Omega} (\bD_k+g_k(|\bF_k|)\bF_k \bF_k^T): \nabla \bw= 0,
\end{aligned}\\
\label{MFweaki}
&\begin{aligned}
\langle \partial_t \bF_k, \bA \rangle &-\int_{\Omega} (\bF_k \otimes \bv_k): \nabla \bA - \int_{\Omega} g_k(|\bF_k|)\nabla \bv_k \bF_k: \bA \\
&+ \frac{1}{2}\int_{\Omega} \left((\bF_k \bF^T_k)^{1+\mu} - (\bF_k \bF^T_k)^{\mu}\right) \bF_k: \bA = 0,
\end{aligned}
\end{align}
and the initial conditions $\bv_0$, $\bF^k_0$ are attained in the following sense:
\begin{equation}
\label{MinvFweaki}
\lim_{t\to 0+} \left(\|\bv_k(t)-\bv_0\|_2 + \|\bF_k(t)-\bF^k_0\|_2\right) = 0.
\end{equation}

To derive the a~priori estimates for the system \eqref{Mvweaki}--\eqref{MFweaki}, we again have to mollify the equation \eqref{MFweaki} and appeal to the Friedrichs commutator lemma (see Section \ref{CFK}, Step 1a, or \cite{BuMaLo22} for details). 
Hence, we set $\bw:=\bv_k$ in \eqref{Mvweaki} and in the mollified form of \eqref{MFweaki} we take the mollification of $\bF_k$ for the test function $\bA$. Then we sum the resulting identities and, after taking the limit in the mollifying parameter, using the fact that $\bv_k\in L^2(0,T;W_{\vec{0},\diver}^{1,2})$, the integration by parts and Korn's and Young's inequalities, we obtain
\begin{equation}\label{1stai}
\begin{split}
\sup_{t\in (0,T)}\left(\|\bv_k(t)\|_2^2 + \|\bF_k(t)\|_2^2 \right) + \int_0^T \|\nabla \bv_k\|_2^2 + \|\bF_k\|_{4+2\mu}^{4+2\mu}\\ \le C(\mu,\Omega)T+\tilde{C}\left(\|\bv_0\|_2^2 + |\bF_0\|_2^2 \right).
\end{split}
\end{equation}
From \eqref{Mvweaki}, \eqref{MFweaki}, \eqref{1stai}, Hölder's and Minkowski's inequalities, standard Sobolev's embeddings and the standard interpolation inequality for the velocity fields, we conclude the estimate
\begin{equation}\label{2stai}
\begin{split}
&\int_0^T \|\partial_t \bv_k\|_{\left(W_{\vec{0},\diver}^{1,\frac{d+2}{2}}\right)^*}^{\frac{d+2}{d}} + \|\partial_t \bF_k\|_{\left((W^{1,d}(\Omega))^{d\times d}\right)^*}^{\frac{4+2\mu}{3+2\mu}}  \\& \qquad \qquad \qquad + \int_0^T \|\bv_k\|_{\frac{2(d+2)}{d}}^{\frac{2(d+2)}{d}} \le C(\mu,T,\Omega,\bF_0,\bv_0).
\end{split}
\end{equation}

Moreover, following step by step the proof of \eqref{3sta}, we deduce that
\begin{equation}\label{blndetFki}
    \sup_{t\in (0,T)}\|\ln \det \bF_k\|_{1}\leq C(\mu, T, \Omega, \bv_0, \bF_0)\quad \implies \quad \det \bF_k>0 \mbox{ a.e. in } Q_T.
\end{equation}

The uniform estimates \eqref{1stai} and \eqref{2stai} together with the definition of $g_k$ and the Aubin--Lions compactness lemma imply the existence of $(\bv, \bF)$ such that, for a  suitable (not relabeled) subsequence of $\{(\bv_k, \bF_k)\}_{k\in \N}$, it holds 
\begin{align}
    \label{vkweaki}\bv_k&\rightharpoonup \bv && \mbox{weakly-* in } L^\infty(0,T;L^2_{\bn, \diver})\cap L^2(0,T; W_{\vec{0},\diver}^{1,2}),\\
    \label{dtvkweaki}\partial_t \bv_k &\rightharpoonup \partial_t \bv && \mbox{weakly in } L^\frac{d+2}{d}\left(0,T; (W_{\vec{0},\diver}^{1,\frac{d+2}{2}})^*\right),\\
    \label{Fkweaki}\bF_k&\rightharpoonup \bF && \mbox{weakly-* in } L^\infty(0,T;(L^2(\Omega))^{d\times d})\cap (L^{4+2\mu}(Q_T))^{d\times d},\\
    \label{dtFkweaki}\partial_t \bF_k &\rightharpoonup \partial_t \bF&& \mbox{weakly in } L^\frac{4+2\mu}{3+2\mu}\left(0,T;((W^{1,d}(\Omega))^{d\times d})^*\right),
\end{align}
and
\begin{align}
    \bv_k&\to \bv && \mbox{strongly in } (L^q(Q_T))^{d\times d}\ \forall q\in \left[1,\frac{2(d+2)}{d}\right),\label{vkstrongi}\\
\label{gkvkstrongi}
    g_k(|\bv_k|)&\to 1 && \mbox{strongly in } L^q(Q_T)\ \forall q\in [1,\infty),\\
    \label{gkFkstrongi}
    g_k(|\bF_k|)&\to 1 && \mbox{strongly in } L^q(Q_T)\ \forall q\in [1,\infty).
\end{align}
Hence, taking the limit $k\to \infty$ in \eqref{Mvweaki}--\eqref{MFweaki}, we obtain, for all $\bw\in W_{\vec{0},\diver}^{1,\frac{d+2}{2}}$, $\bA \in (W^{1, d}(\Omega))^{d\times d}$ and a.a. $t\in (0,T)$, that
\begin{align} \label{2wherei}&\langle\partial_t \bv, \bw\rangle - \int_{\Omega} (\bv \otimes \bv): \nabla\bw + \int_{\Omega} \bD: \nabla \bw + \int_{\Omega} \overline{\bF \bF^T}: \nabla \bw = 0, \\ &
\label{Fwithoutcompi}
\begin{aligned}
\langle\partial_t \bF, \bA\rangle &- \int_{\Omega}(\bF \otimes\bv): \nabla\bA - \int_{\Omega} \overline{(\nabla \bv) \bF}:\bA 
\\ &+ \frac{1}{2}\int_{\Omega}\left(\overline{(\bF \bF^T)^{1+\mu}\bF}  -\overline{(\bF\bF^T)^{\mu}\bF}\right):\bA = 0.
\end{aligned}
\end{align}
The proofs of the properties $\bv\in \mathcal{C}_{weak}([0,T];L^2_{\bn, \diver})$, $\bF\in \mathcal{C}([0,T]; (L^2(\Omega))^{d\times d})$ and of the attainment initial conditions \eqref{invFweaki} can be read directly from the proofs of the corresponding properties for the Giesekus model, hence we do not present them here. It remains to prove that $\ln \det \bF\in L^\infty(0,T;L^1(\Omega))$, hence $\det \bF>0$ a.e. in $Q_T$, and to identify that $\overline{\bF \bF^T}=\bF \bF^T$ in \eqref{2wherei} and $\overline{\nabla \bv \bF}=\nabla \bv \bF$, $\overline{(\bF \bF^T)^{1+\mu} \bF}=(\bF \bF^T)^{1+\mu} \bF$ and $\overline{(\bF \bF^T)^{\mu} \bF}=(\bF \bF^T)^{\mu} \bF$ in \eqref{Fwithoutcompi}. For the identification of the weak limits it suffices, due to the convergence results \eqref{vkweaki}--\eqref{gkFkstrongi}, to prove the compactness of $\{\bF_k\}_{k\in \N}$ in $(L^2(Q_T))^{d\times d}$. The property $\ln \det \bF\in L^\infty(0,T;L^1(\Omega))$, hence $\det \bF>0$ a.e. in $Q_T$, is a direct consequence of the compactness of $\{\bF_k\}_{k\in \N}$ in $(L^2(Q_T))^{d\times d}$, \eqref{blndetFki} and the Fatou lemma, see Section \ref{limitk2} for details.

 \subsubsection{Compactness of \texorpdfstring{$\{\bF_k\}_{k\in \N}$ in $(L^2(Q_T))^{d\times d}$}{ed}}
We follow the proof of the compactness of $\{\bF_k\}_{k\in \N}$ in $(L^2(Q_T))^{3\times 3}$ for the Giesekus model. As clarified in Sect.~\ref{CFK}, it suffices to show, for all nonnegative $\varphi\in \mathcal{C}_c^\infty((-\infty,T)\times \Omega)$, the following inequality:
\begin{equation}\label{prepGroni}
    -\int_{Q_T}(\overline{|\bF|^2}-|\bF|^2)\partial_t \varphi-\int_{Q_T} (\overline{|\bF|^2}-|\bF|^2) \bv\cdot \nabla \varphi\leq \int_{Q_T} L(\overline{|\bF|^2}-|\bF|^2)\varphi,
\end{equation}
where $L$ is an $L^2(Q_T)$ function specified below. 

Following Step 1a and Step 1b in Sect.~\ref{CFK} we first arrive at the inequality
\begin{equation}
    \label{step1compFsigmai}
    \begin{split}
    -\int_{Q_T} \left(\overline{|\bF|^2}-|\bF|^2\right)\partial_t \varphi - \int_{Q_T} \left(\overline{|\bF|^2}-|\bF|^2\right) \bv\cdot \nabla \varphi 
     \leq -l_{1}+\tilde{l}_1+2l_{2},
\end{split}
\end{equation}
where
\begin{align}
\label{lbi}
 l_{1}&:= \lim_{m\to \infty}\lim_{k\to \infty} \int_{Q_T} \left(\mbox{tr }(\bF_k \bF_k^T)^{2+\mu}g_m(|\bF_k|^2)-(\bF_k\bF_k^T)^{1+\mu}\bF_k:\bF\right)\varphi,
    \\
    \label{lbti}\tilde{l}_{1}&:= \lim_{m\to \infty}\lim_{k\to \infty} \int_{Q_T} \left(\mbox{tr }(\bF_k \bF_k^T)^{1+\mu}g_m(|\bF_k|^2)-(\bF_k\bF_k^T)^{\mu}\bF_k:\bF\right)\varphi,\\
\label{lai}
    l_{2}&:= \lim_{m\to \infty}\lim_{k\to \infty} \int_{Q_T} \left(\nabla \bv_k: g_m(|\bF_k|^2) g_k(|\bF_k|)\bF_k \bF_k^T- \nabla\bv_k\bF_k:\bF\right)\varphi.
\end{align}
Employing the definition of $g_m$ and the weak convergence \eqref{Fkweaki}, we get, analogously as in  Section \ref{CFK}, Step 2a, that
\begin{align}\label{larew}
l_1&= \lim_{m\to \infty} \lim_{k\to\infty} \int_{Q_T} \left((\bF_k \bF_k^T)^{1+\mu} \bF_k - (\bF \bF^T)^{1+\mu}\bF\right):(\bF_k-\bF)g_m(|\bF_k|^2)\varphi,
\\
\label{lbrew}
\tilde{l}_1&= \lim_{m\to \infty} \lim_{k\to\infty} \int_{Q_T} \left((\bF_k \bF_k^T)^{\mu} \bF_k - (\bF \bF^T)^{\mu}\bF\right):(\bF_k-\bF)g_m(|\bF_k|^2)\varphi.
\end{align}
By Proposition \ref{mon} we see that both $l_1$, $\tilde{l}_1$ are nonnegative and that
\begin{equation}
    -l_1+\tilde{l}_1\leq c(\mu) \lim_{m\to \infty} \lim_{k\to \infty} \int_{Q_T}|\bF_k-\bF|^2 g_m(|\bF_k|^2)\varphi,
\end{equation}
which, due to \eqref{Fkweaki} and the definition of $g_m$, can be written as
\begin{equation}\label{monik}
    -l_1+\tilde{l}_1\leq c(\mu) \int_{Q_T}\left(\overline{|\bF|^2}-|\bF|^2\right) \varphi.
\end{equation}
Repeating the procedure from Section \ref{CFK}, Step 2b, we get
\begin{equation}\label{step2bi}
    l_2\leq \int_{Q_T}(|\bD|+|\bF|^2)\left(\overline{|\bF|^2}-|\bF|^2\right) \varphi.
\end{equation}
Summing \eqref{monik} and \eqref{step2bi} and using \eqref{step1compFsigmai}, we obtain \eqref{prepGroni} with $L:= c(\mu)+2(|\bD|+|\bF|^2)\in L^2(Q_T)$. This completes the proof of the compactness of $\{\bF_k\}_{k\in \N}$ in $(L^2(Q_T))^{d\times d}$.
\qed

\bibliographystyle{plainnat}

\bibliography{name}

\end{document}